# Consensus of a class of nonlinear fractional-order multi-agent systems via dynamic output feedback controller

Elyar Zavvari*, Pouya Badri, Mahdi Sojoodi

Advanced Control Systems Laboratory, School of Electrical and Computer Engineering, Tarbiat Modares University, Tehran, Iran.

*Abstract—This paper addresses the consensus of a class of nonlinear fractional-order multi-agent systems (FOMAS) with positive real uncertainty. First a fractional non-fragile dynamic output feedback controller is put forward via the output measurements of neighboring agents, then appropriate state transformation reduced the consensus problem to a stability one. A sufficient condition based on direct Lyapunov approach, for the robust asymptotic stability of the transformed system and subsequently for the consensus of the main system is presented. Additionally, utilizing S-procedure and Schur complement, the systematic stabilization design algorithm is proposed for fractional-order system with and without nonlinear terms. The results are formulated as an optimization problem with linear matrix inequality constraints. Simulation results are given to verify the effectiveness of the theoretical results.*

## I. INTRODUCTION

Appearance of fractional calculus emerged the idea of modeling systems via non-integer order differential operators. Most of real systems mainly have fractional behavior, so it could be worthwhile to describe them with fractional operators. Fractional calculus developed new mathematical tools better describing real-world systems, in comparison with traditional integer-order derivative equations [1], [2]. A basic issue in control theory to develop solutions for control objectives is to have an accurate model of the systems [3]–[6]. The absence of fractional-order differential equations was the main reason for using integer-order models in control theory. The emergence of methods for approximation of fractional derivative and integral paved the way for using fractional calculus in wide areas of control theory. Some examples of fractional systems include viscoelastic polymers [7], biomedical applications [8], semi-infinite transmission lines with losses [9], dielectric polarization [10].

Modeling and study of multi-agent systems have attracted tremendous attention in recent years [11]–[13]. This is partly due to their potential applications in many areas, including control theory, mathematics, biology, physics, computer science, and robotics. Consensus is the concept of reaching an agreement considering the states of all agents [14], [15] and plays an important role in multi-agent systems. Examples include consensus of a class of nonlinear systems with dynamic output feedback [14], flocking [16], formation control [17], cooperative control [18], distributed sensor networks [19], synchronization between the motors, and so on.

---

* Corresponding Author, Address: P.O. Box 14115-…, Tehran, Iran, E-mail address: elyar.zavary@gmail.com, Tel./Fax: +98 21 8288-3902.

The problem of robust consensus of fractional-order linear multi-agent systems via static feedback was studied in [20], furthermore, [21] investigates the distributed containment control of fractional-order uncertain multi-agent systems. Control and synchronization of a class of uncertain fractional-order chaotic systems via adaptive backstepping Control is studied in [22]. Consensus control of fractional-order systems based on delayed state fractional-order derivative is investigated in [15]. Then, the static output feedback controller is utilized to stabilize the transformed system. It is worth mentioning that controllers, designed based on dynamic feedback, are always preferable to the static ones because of their more effective control performances, moreover, the dynamic controller brings about more degree of freedom in achieving control objectives, in comparison with the static controller [23]. In addition, most of mentioned works use state feedback controller and this kind of controllers require all states. On the other hand, in some cases states are inaccessible because of costly implementation or some physical constraints. High-order controllers obtained by most of the controller design methods have expensive implementation procedures, undesirable reliability, high fragility, and numerous maintenance difficulties. Since controller order reduction techniques may deteriorate the closed-loop efficiency, designing directly a low-, fixed-order controller for a system can be helpful [24]–[26].

The fractional-order multi-agent system defined in this paper is affected by nonlinear uncertainty. The main purpose of the paper is to derive some LMI constraints that can be easily checked by existing solvers and parsers to obtain the consensus controller parameters. It does not seem straightforward to derive linear stabilizing constraints for a linear system which is affected by nonlinear uncertainties. However, using helpful lemmas which were introduced in Preliminaries and problem formulation section, and also utilizing LMI techniques which were used in the proof procedures of the theorems, nonlinear uncertain parameters were involved in the achieved LMI constraints. Then, by solving the mentioned LMIs, the unknown parameters of the controller can be achieved to establish consensus of the uncertain multi-agent system. Therefore, stabilizing control systems plays a crucial role in control engineering is an undeniable fact.

Moreover, uncertainty is in real systems' nature. As a result, a dynamic robust controller covers a wide range of systems and its advantages overtake static controller ones [11], [20], [24]. Besides that, designing a dynamic robust controller leads to more unknown parameters in comparison with a static controller and makes controller design procedure more difficult due to more complex constraints which must be solved by the solvers. In this paper, using proper lemmas and theorems, LMI techniques, and suitable solvers and parsers the difficulty of designing such controllers has been overcome.

Numerous complex phenomena in many real-world systems and their behavior can be well described by nonlinear dynamics. Existence of nonlinearity in most of systems, motivates researchers to investigate control and stability methods on the system with nonlinear behavior. As a result, studying nonlinear FOMASs cooperative behaviour is useful. The consensus of nonlinear FOMASs was studied in [27]–[32]. A Fractional-order complex networks with Lipschitz-type nonlinear nodes and directed communication topology had been investigated in [31], this paper mainly focused on pinning synchronization problem. The consensus of all agents with nonlinear terms containing internal and coupling delays is guaranteed with a heterogeneous impulsive control method in [27]. By means of Mittag-Leffler function, the Laplace transform, and the inequality technique [32] had studied the leader-following consensus of nonlinear FOMASs. Robust consensus tracking problem of FOMASs in the presence of heterogeneous unknown nonlinearities and external disturbances and output feedback consensus of nonlinear FOMASs with directed topologies were investigated in [30] and [28] respectively.

However, all mentioned works deal with cooperative behavior of nonlinear FOMASs with a special system and input matrices. Designing non-fragile dynamic output feedback controller for the consensus of nonlinear fractional-order multi-agent systems (FOMAS) with the general linear dynamics for the linear part, which is more general, is addressed in this work. The main contributions of this paper are summarized as follows.

- Unlike many existing implementations [1], [2], [15], [33], [34] which are limited to linear FOMAS, while our proposed framework considers a more general class of agents with general nonlinear fractional order dynamics.
- Given that only gain uncertainty can be obtained by polytopic uncertainty, norm-bounded uncertainty, and interval uncertainty descriptions and consequently conservative results will be captured. Applying the positivity theorem and modeling the uncertainty with a real positive model is a way for accounting phase information. As a result, this model of system uncertainty description not only adds a certain degree of robustness in our proposed framework but also prevents conservative results.
- Compared to the commonly used state feedback control protocols for consensus in FOMAS [11], [12], [35], the implementation in this article is based on the dynamic output feedback protocol with its well-known benefits.
- Additionally, our design method allows a priory degree of uncertainty in the controller gains so that the implementation reflects resilience to uncertainties from the controller perspective

***Notations:*** *In this paper,* $A \otimes B$ *denotes the kronecker product of matrices A and B, and the symmetric of matrix M will be shown by* $sym(.)$*, which is defined by* $sym(M) = M^T + M$*, and also* $\uparrow$ *is the symbol of pseudo inverse of matrix.*

## II. Problem Formulation and Preliminaries

In this section, some basic notations and definitions related to fractional-order systems are given. Some concepts and lemmas about graph theory are presented as well. $\boldsymbol{I}_n$ denotes the $n \times n$ identity matrix, $\mathbf{1}_n$ and $\mathbf{0}_n$ indicate $n \times 1$ column vectors with all elements to be ones and zeros, respectively. Moreover, $\boldsymbol{J}_n$ is a matrix of $n \times n$ dimension, with all elements to be ones. $A \circ B$ is a matrix with elements defined by

$$(A \circ B)_{ij} = (A)_{ij}.(B)_{ij}.$$

A team of $N$ ($N > 1$) fractional-order networked nonlinear agent systems are considered. The dynamic of the $i$th agent is described as

$$\begin{aligned} D^q \boldsymbol{x}_i(t) &= \left(\widetilde{\boldsymbol{\mathcal{A}}} + \boldsymbol{\Delta}\widetilde{\boldsymbol{\mathcal{A}}}_i\right)\boldsymbol{x}_i(t) + \widetilde{\boldsymbol{\mathcal{B}}}_i \boldsymbol{u}_i(t) + \widetilde{\boldsymbol{\phi}}\left(\boldsymbol{x}_i(t), \boldsymbol{u}_i(t)\right) \\ \boldsymbol{y}_i(t) &= \widetilde{\boldsymbol{\mathcal{C}}}\boldsymbol{x}_i(t), \qquad i = 1, \dots, N \end{aligned} \tag{1}$$

with initial condition

$$\boldsymbol{x}_i(0) = \boldsymbol{x}_{i0} \tag{2}$$

where $\boldsymbol{x}_i(t) \in \mathbb{R}^n, \boldsymbol{u}_i(t) \in \Re^m, \boldsymbol{y}_i(t) \in \mathbb{R}^p$ are pseudo state, input, measured output, respectively. $\widetilde{\boldsymbol{\mathcal{A}}} \in \mathbb{R}^{n\times n}, \widetilde{\boldsymbol{\mathcal{B}}}_i \in \mathbb{R}^{n\times m}, \widetilde{\boldsymbol{\mathcal{C}}} \in \Re^{p\times n}$ are known constant matrices, and $\widetilde{\boldsymbol{\phi}}(\cdot) : [\mathbb{R}^n \times \mathbb{R}^m] \to \mathbb{R}^n,$ is nonlinear function. $\boldsymbol{\Delta}\widetilde{\boldsymbol{\mathcal{A}}}_i \in \mathbb{R}^{n\times n}$ is a time-invariant matrix, with parametric uncertainty. $q$ is the fractional derivative order, there are several definitions for fractional-order derivative, among them Grünwald-Letnikov, Riemann-Liouville, and Caputo are most commonly referred. However, since Caputo definitions initial condition is similar to integer orders one as a physical aspect, Caputo definition is used in this with the following definition

$${}_{a}^{C}D_t^q = \frac{1}{\Gamma(\bar{n}-\alpha)} \int_a^t (t-\tau)^{\bar{n}-a-1} \left(\frac{d}{d\tau}\right)^{\bar{n}} f(\tau)d\tau,$$

where $\boldsymbol{\Gamma}(\cdot)$ is Gamma function defined by $\boldsymbol{\Gamma}(\epsilon) = \int_0^\infty e^{-t}t^{\epsilon-1}dt$ and $\bar{n}$ is the smallest integer that is equal or greater than $q$.

***Lemma 1.*** *([36]) Let* $\boldsymbol{\mathcal{A}} \in \mathbb{R}^{n\times n}$, $0 < q < 1$ *and* $\theta = q\pi/2$. *The fractional-order system* $D^q \boldsymbol{x}(t) = \boldsymbol{\mathcal{A}}\boldsymbol{x}(t)$ *is asymptotically stable if and only if there exist two real symmetric positive definite matrices* $X_{k1} \in \mathbb{R}^{n\times n}$, $k = 1,2$, *and two skew-symmetric matrices* $X_{k2} \in \mathcal{R}^{n\times n}$, $k = 1,2$, *such that*

$$\begin{aligned} &\sum_{i=1}^{2}\sum_{j=1}^{2} Sym\{\Theta_{ij} \otimes (\boldsymbol{\mathcal{A}}X_{ij})\} < 0 \\ &\begin{bmatrix} X_{11} & X_{12} \\ -X_{12} & X_{11} \end{bmatrix} > 0 \quad \begin{bmatrix} X_{21} & X_{22} \\ -X_{22} & X_{21} \end{bmatrix} > 0 \end{aligned} \tag{3}$$

*where*

$$\begin{aligned} &\Theta_{11} = \begin{bmatrix} \sin\theta & -\cos\theta \\ \cos\theta & \sin\theta \end{bmatrix}, \quad \Theta_{12} = \begin{bmatrix} \cos\theta & \sin\theta \\ -\sin\theta & \cos\theta \end{bmatrix} \\ &\Theta_{21} = \begin{bmatrix} \sin\theta & \cos\theta \\ -\cos\theta & \sin\theta \end{bmatrix}, \quad \Theta_{22} = \begin{bmatrix} -\cos\theta & \sin\theta \\ -\sin\theta & -\cos\theta \end{bmatrix}. \end{aligned} \tag{4}$$

***Lemma 2.*** *[37] Let* $f: \mathbb{R}_\epsilon \to \mathbb{R}^n$ *be piecewise continuous respect to* $t$, *where* $\mathbb{R}_\epsilon = \{(t, \boldsymbol{x}): 0 \le t \le a \text{ and } \|\boldsymbol{x} - \boldsymbol{x_0}\| \le b\}$, $f = [f_1, \dots, f_n]^T$, $\boldsymbol{x} \in \mathbb{R}^n$ *and* $\|f(t, \boldsymbol{x})\| \le M$ *on* $\mathbb{R}_\epsilon$. *Then, there exists at least one solution for the system of fractional differential equations given by*

$$D^q \boldsymbol{x}(t) = \boldsymbol{f}(t, \boldsymbol{x}(t)) \tag{5}$$

*with the initial condition*

$$\boldsymbol{x}(0) = \boldsymbol{x}_0 \tag{6}$$

*on* $0 \le t \le \beta$ *where* $\beta = min\left(a, \left[(b/M)\Gamma(q+1)^{1/q}\right]\right), 0 < q < 1.$

***Lemma 3.*** *[37] Consider initial fractional problem (5) and (6) with* $0 < q < 1$ *and assume that conditions of Lemma 2 hold. Let*

$$g(v, \boldsymbol{x}_*(v)) = f\left(t - \left(t^q - v\Gamma(q+1)\right)^{1/q}, \boldsymbol{x}\left(t - v\Gamma(q+1)\right)^{1/q}\right)$$

*then* $\boldsymbol{x}(t)$, *is given by*

$$\boldsymbol{x}(t) = \boldsymbol{x}_*(t^q/\Gamma(\boldsymbol{q}+\boldsymbol{1})),$$

*where* $\boldsymbol{x}_*(v)$ *can be obtained by solving the following integer order differential equation*

$$\begin{aligned} &\frac{d\boldsymbol{x}_*(v)}{dv} = g(v, \boldsymbol{x}_*(v)) \\ &\boldsymbol{x}_*(0) = \boldsymbol{x}_0. \end{aligned} \tag{7}$$

System matrices $\widetilde{\boldsymbol{\mathcal{A}}}$, $\widetilde{\boldsymbol{\mathcal{B}}}_i$, $\widetilde{\boldsymbol{\mathcal{C}}}$, nonlinear function $\boldsymbol{\phi}(\cdot)$ and uncertainty matrix $\boldsymbol{\Delta}\widetilde{\boldsymbol{\mathcal{A}}}_i$ are assumed to satisfy the following assumptions.

***Remark 1.*** *An integro-differential Laplacian operator* $s^{\pm\alpha}$, *has the frequency response of the form* $(j\omega)^{\pm\alpha} = \omega^{\pm\alpha}\left[\cos\left(\frac{\alpha\pi}{2}\right) \pm j\sin(\alpha\pi/2)\right]$. *Analytical inverse Laplace transform of a fractional-order transfer function does not exist* [38]. *Consequently, there is no way to estimate the exact time response of a fractional-order transfer function but approximations. There are some popular approximations to deal with fractional order behavior of systems, like least square method, Oustaloup's, Continued Fractional Expansion, Carlson's and Matsuda's methods. Furthermore, by the means of system's impulse and the step responses, the fractional order system can be modeled by an integer order system. There is a limitation using this approach that the approximated integer order system is valid only for a specific frequency band. In this paper Grünwald–Letnikov method is used to approximate the fractional derivatives. Grünwald–Letnikov method is another approach which is widely used in different research areas to approximate fractional order derivative* *and for numerical method's implementation is very effective* [5], [13], [39]. *As a result, we have used Grünwald–Letnikov approximation approach for numerical solutions in this paper.*

***Assumption 1.*** *The pairs of* $(\widetilde{\mathcal{A}}, \widetilde{\mathcal{B}}_i)$ *and* $(\widetilde{\mathcal{A}}, \widetilde{\mathcal{C}})$ *are controllable and observable, respectively.*

***Assumption 2.*** $\Delta\widetilde{\mathcal{A}}_i$ *is time-invariant matrix of the following form:*

$$\Delta\widetilde{\mathcal{A}}_i = \widetilde{\mathcal{R}}\delta_i(\sigma)\widetilde{\mathcal{N}}, \tag{8}$$

$$\delta_i(\sigma) = \mathcal{Z}_i(\sigma)[\boldsymbol{I} + \mathcal{J}\mathcal{Z}_i(\sigma)]^{-1}, \tag{9}$$

$$Sym\{\mathcal{J}\} > 0, \tag{10}$$

*where* $\widetilde{\mathcal{R}} \in \mathbb{R}^{n\times m_0}$, $\widetilde{\mathcal{N}} \in \mathbb{R}^{m_0\times n}$, *and* $\mathcal{J} \in \mathbb{R}^{m_0\times m_0}$ *are real known matrices. The uncertain matrix* $\mathcal{Z}_i(\sigma) \in \mathbb{R}^{m_0\times m_0}$ *satisfies*

$$Sym\{\mathcal{Z}_i(\sigma)\} \geq 0, \tag{11}$$

*where* $\sigma \in \Omega$, *with* $\Omega$ *being a compact set.*

***Remark 2.*** *Condition (10) guarantees that* $I + JZ_i(\sigma)$ *is invertible for all* $Z_i(\sigma)$ *satisfying (11), therefore* $\delta_i(\sigma)$ *in (8) is well defined ([40]).*

***Assumption 3.*** *Nonlinear function* $\widetilde{\boldsymbol{\phi}}(\boldsymbol{x}_i(t), \boldsymbol{u}_i(t))$ *is Lipschitz on* $\boldsymbol{x}_i(t)$ *with Lipschitz constant* $\xi_1$

$$\|\widetilde{\boldsymbol{\phi}}(\boldsymbol{x}_1(t), \boldsymbol{u}_1(t)) - \widetilde{\boldsymbol{\phi}}(\boldsymbol{x}_2(t), \boldsymbol{u}_2(t))\| < \xi_1\|\boldsymbol{x}_1(t) - \boldsymbol{x}_2(t)\| \tag{12}$$

*for all* $\boldsymbol{x}_1(t), \boldsymbol{x}_2(t) \in \mathbb{R}^n$ *and*

$$\widetilde{\boldsymbol{\phi}}(0,0) = 0 \tag{13}$$

The graph $\boldsymbol{\mathcal{G}} = (\boldsymbol{\mathcal{V}}, \varepsilon, \boldsymbol{A})$ describes the information exchanging among agents, where $\boldsymbol{\mathcal{V}} = \{v_1, \dots, v_p\}$ is the vertex set, $\varepsilon \subseteq \boldsymbol{\mathcal{V}} \times \boldsymbol{\mathcal{V}}$ is the edge set, and $\boldsymbol{A} = (a_{ij})_{p\times p}$ is a nonsymmetrical set. A nonzero $\varepsilon_{ij} \in \varepsilon \subseteq \boldsymbol{\mathcal{V}} \times \boldsymbol{\mathcal{V}}$ indicates that agent $j$ receives information from agent $i$ which leads to a corresponding nonzero $a_{ij} \in \boldsymbol{A} = (a_{ij})_{p\times p}$ and $a_{ij} = 0$ otherwise. Furthermore, $a_{ii}$ is supposed to be zero for all $i \in \{1, \dots, n\}$ and $\boldsymbol{\mathcal{N}}_i = \{v_j | \varepsilon_{ij} \in \varepsilon\}$ is the set of neighbors of agent $i$ [41], as well as $\boldsymbol{L}_{n\times n}$, is the Laplacian matrix of graph $\boldsymbol{\mathcal{G}}$, which is defined as follows

$$\boldsymbol{L} = (l_{ij}), \quad l_{ij} = \begin{cases} \sum_{p\in\boldsymbol{\mathcal{N}}_i} a_{ip}, & i = j \\ -a_{ij}, & i \neq j \end{cases}. \tag{14}$$

***Lemma 4.*** *([42]) Let* $\Omega = \{\delta \in \mathbb{R}^{m_0\times m_0} | \delta \text{ is subjected to } (8) - (10)\}$. *Then*
$\Omega = \{\delta \in \mathbb{R}^{m_0\times m_0} | \det(I - \delta J) \neq 0 \text{ and } \delta Sym\{J\}\delta^T \leq Sym\{\delta\}\}$.

## III. MAIN RESULT

In this work, we will study the consensus problem for FOMASs composed of the system defined in *(1)*. The multi-agent system *(1)* can be represented in the following augmented form

$$\begin{aligned} &D^\alpha \boldsymbol{x}(t) = \overline{\boldsymbol{\mathcal{A}}}_{\Delta,N}\boldsymbol{x}(t) + \boldsymbol{\mathcal{B}}\boldsymbol{u}(t) + \boldsymbol{\phi}(\boldsymbol{x}(t), \boldsymbol{u}(t)), \\ &\boldsymbol{y}(t) = \boldsymbol{\mathcal{C}}\boldsymbol{x}(t), \\ &\boldsymbol{x}(0) = \boldsymbol{x}_0 \end{aligned} \tag{15}$$

where $\boldsymbol{x}(t) = [\boldsymbol{x}_1(t)^T, \dots, \boldsymbol{x}_N(t)^T]^T \in \mathbb{R}^{N.n}$, $\boldsymbol{u}(t) = [\boldsymbol{u}_1(t)^T, \dots, \boldsymbol{u}_N(t)^T]^T \in \mathbb{R}^{N.m}$, and $\boldsymbol{y}(t) = [\boldsymbol{y}_1(t)^T, \dots, \boldsymbol{y}_N(t)^T]^T \in \mathbb{R}^{N.p}$ are general $x$pseudo state, input, and output vectors respectively, and also $\overline{\boldsymbol{\mathcal{A}}}_{\Delta,N} = \boldsymbol{\mathcal{A}}_N +$ $\boldsymbol{\Delta\mathcal{A}_N}$ where $\boldsymbol{\mathcal{A}}_N = \boldsymbol{I}_N \otimes \widetilde{\boldsymbol{\mathcal{A}}}$ and $\boldsymbol{\Delta\mathcal{A}_N} = diag(\boldsymbol{\Delta}\widetilde{\boldsymbol{\mathcal{A}}}_i)$, $\boldsymbol{\mathcal{B}} = diag(\widetilde{\boldsymbol{\mathcal{B}}}_1, \dots, \widetilde{\boldsymbol{\mathcal{B}}}_N) \in \mathbb{R}^{N.n\times N.m}$, and $\boldsymbol{\mathcal{C}}_N = \boldsymbol{I}_N \otimes \widetilde{\boldsymbol{\mathcal{C}}}$ are constant known matrices. $\boldsymbol{\phi}(\boldsymbol{x}(t), \boldsymbol{u}(t)) = [\widetilde{\boldsymbol{\phi}}^T(\boldsymbol{x}_1(t), \boldsymbol{u}_1(t)), \dots, \widetilde{\boldsymbol{\phi}}^T(\boldsymbol{x}_N(t), \boldsymbol{u}_N(t))]^T$ is augmented form of system nonlinear term. Our aim is to define a fractional-order decentralized dynamic output feedback controller that guarantees consensus of agents.

***Definition 1.*** *System (15) achieves consensus asymptotically if the following condition holds for any* $\boldsymbol{x}_i(0) = \boldsymbol{x}_{i0}$

$$\lim_{t\to\infty} |\boldsymbol{x}_i(t) - \boldsymbol{x}_j(t)| = 0, \quad i, j = 1, \dots, N, i \neq j. \tag{16}$$

***Definition 2.*** *Consensus error of ith agent of system (1) is defined as follows*

$$\mathbf{e}_i(t) = \sum_{\boldsymbol{j=1,\dots,N}} \|\boldsymbol{x}_i(t) - \boldsymbol{x}_j(t)\| \ (j \neq i) \tag{17}$$

In order to solve the consensus problem for system *(15)*, we use the following non-fragile control protocol

$$\begin{aligned} &D^q \boldsymbol{x}_{ci}(t) = \overline{\boldsymbol{\mathcal{A}}}_{ci}\, \boldsymbol{x}_{ci}(t) + \overline{\boldsymbol{\mathcal{B}}}_{ci}\big(l_{ii}\boldsymbol{y}_i(t) + \textstyle\sum_{o\in\boldsymbol{\mathcal{N}}_i} l_{io}\boldsymbol{y}_o(t)\big), \\ &\boldsymbol{u}_i = \overline{\boldsymbol{\mathcal{C}}}_{ci}\boldsymbol{x}_{ci} + \overline{\boldsymbol{\mathcal{D}}}_{ci}\big(l_{ii}\boldsymbol{y}_i(t) + \textstyle\sum_{o\in\boldsymbol{\mathcal{N}}_i} l_{io}\boldsymbol{y}_o(t)\big), \qquad i = 1, \dots, N. \\ &\boldsymbol{x}_{ci}(0) = \boldsymbol{x}_{ci0} \end{aligned} \tag{18}$$

where $\boldsymbol{x}_{ci} \in \mathbb{R}^{n_c}$ is controller pseudo state in which $n_c$ is the controller order and $\boldsymbol{\mathcal{N}}_i$ denotes the neighbors of $i$th agent.

***Remark 3.*** *In state feedback control scheme, all individual states of the system are needed to be measured and used in feedback line. However, in some practical situations measuring all states is impossible or may sound difficult due to economic issues or physical limitations [39], [43], [44]. In these cases, using output feedback control could be effective since there is no need to measure all individual states of the system and only by measuring outputs of the system the control action is done.*

***Remark 4.*** *In terms of the feasibility of the solution for control design, it is shown in [45] that the domain of feasibility increases by having additional parameters in solving an inequality. In another word, when we have more parameters in solving a linear matrix inequality our degree of freedom increases. Accordingly, the proposed dynamic output feedback controller guarantees the consensus of the agents more effectively and solver has more degree of freedom in solving the inequality in comparison to a static feedback controller.*

The system matrices $\overline{\boldsymbol{\mathcal{A}}}_{ci}$, $\overline{\boldsymbol{\mathcal{B}}}_{ci}$, $\overline{\boldsymbol{\mathcal{C}}}_{ci}$ and $\overline{\boldsymbol{\mathcal{D}}}_{ci}$ have the following admissible time variant uncertainty:

$$\begin{aligned} &\overline{\boldsymbol{\mathcal{A}}}_{ci} = \boldsymbol{\mathcal{A}}_{ci} + \boldsymbol{D}_{\boldsymbol{\mathcal{A}}_{ci}}\boldsymbol{F}_{\boldsymbol{\mathcal{A}}_{ci}}(t)\boldsymbol{E}_{\boldsymbol{\mathcal{A}}_{ci}} \\ &\overline{\boldsymbol{\mathcal{B}}}_{ci} = \boldsymbol{\mathcal{B}}_{ci} + \boldsymbol{D}_{\boldsymbol{\mathcal{B}}_{ci}}\boldsymbol{F}_{\boldsymbol{\mathcal{B}}_{ci}}(t)\boldsymbol{E}_{\boldsymbol{\mathcal{B}}_{ci}}, \\ &\overline{\boldsymbol{\mathcal{C}}}_{ci} = \boldsymbol{\mathcal{C}}_{ci} + \boldsymbol{D}_{\boldsymbol{\mathcal{C}}_{ci}}\boldsymbol{F}_{\boldsymbol{\mathcal{C}}_{ci}}(t)\boldsymbol{E}_{\boldsymbol{\mathcal{C}}_{ci}}, \\ &\overline{\boldsymbol{\mathcal{D}}}_{ci} = \boldsymbol{\mathcal{D}}_{ci} + \boldsymbol{D}_{\boldsymbol{\mathcal{D}}_{ci}}\boldsymbol{F}_{\boldsymbol{\mathcal{D}}_{ci}}(t)\boldsymbol{E}_{\boldsymbol{\mathcal{D}}_{ci}}, \qquad i = 1, \dots, N. \end{aligned}$$

where $\boldsymbol{D}_{\boldsymbol{\mathcal{A}}_{ci}}$, $\boldsymbol{E}_{\boldsymbol{\mathcal{A}}_{ci}}$, $\boldsymbol{D}_{\boldsymbol{\mathcal{B}}_{ci}}$, $\boldsymbol{E}_{\boldsymbol{\mathcal{B}}_{ci}}$, $\boldsymbol{D}_{\boldsymbol{\mathcal{C}}_{ci}}$, $\boldsymbol{E}_{\boldsymbol{\mathcal{C}}_{ci}}$, $\boldsymbol{D}_{\boldsymbol{\mathcal{D}}_{ci}}$ and $\boldsymbol{E}_{\boldsymbol{\mathcal{D}}_{ci}}$ are known constant matrices, and $\boldsymbol{F}_{\boldsymbol{\mathcal{A}}_{ci}}(t)$, $\boldsymbol{F}_{\boldsymbol{\mathcal{B}}_{ci}}(t)$, $\boldsymbol{F}_{\boldsymbol{\mathcal{C}}_{ci}}(t)$ and $\boldsymbol{F}_{\boldsymbol{\mathcal{D}}_{ci}}(t)$ are unknown matrices Lebesgue measurable elements satisfying

$$\begin{aligned} &\boldsymbol{F}^T_{\boldsymbol{\mathcal{A}}_{ci}}(t)\boldsymbol{F}_{\boldsymbol{\mathcal{A}}_{ci}}(t) < \boldsymbol{I} \\ &\boldsymbol{F}^T_{\boldsymbol{\mathcal{B}}_{ci}}(t)\boldsymbol{F}_{\boldsymbol{\mathcal{B}}_{ci}}(t) < \boldsymbol{I} \\ &\boldsymbol{F}^T_{\boldsymbol{\mathcal{C}}_{ci}}(t)\boldsymbol{F}_{\boldsymbol{\mathcal{C}}_{ci}}(t) < \boldsymbol{I} \\ &\boldsymbol{F}^T_{\boldsymbol{\mathcal{D}}_{ci}}(t)\boldsymbol{F}_{\boldsymbol{\mathcal{D}}_{ci}}(t) < \boldsymbol{I} \qquad i = 1, \dots, N. \end{aligned} \tag{19}$$

Defining $\boldsymbol{x}_c = [\boldsymbol{x}_{c1}^T \quad \dots \quad \boldsymbol{x}_{cN}^T]^T \in \mathbb{R}^{N.n_c}$ as controller pseudo state vector, $\boldsymbol{\mathcal{A}}_c = diag(\boldsymbol{\mathcal{A}}_{c1}, \dots, \boldsymbol{\mathcal{A}}_{cN}) \in \mathbb{R}^{N.n_c \times N.n_c}$, $\boldsymbol{D}_{\mathcal{A}_c} = diag\left(\boldsymbol{D}_{\mathcal{A}_{c1}}, \dots, \boldsymbol{D}_{\mathcal{A}_{cN}}\right)$ $\boldsymbol{\mathcal{B}}_c = diag(\boldsymbol{\mathcal{B}}_{c1}, \dots, \boldsymbol{\mathcal{B}}_{cN}) \in \mathbb{R}^{N.n_c \times N.p}$, $\boldsymbol{\mathcal{C}}_c = diag(\boldsymbol{\mathcal{C}}_{c1}, \dots, \boldsymbol{\mathcal{C}}_{cN}) \in \mathbb{R}^{N.m \times N.n_c}$, $\boldsymbol{\mathcal{D}}_c = diag(\boldsymbol{\mathcal{D}}_{c1}, \dots, \boldsymbol{\mathcal{D}}_{cN}) \in \mathbb{R}^{N.m \times N.p}$ as controller matrices, and

$$
\begin{aligned}
&\boldsymbol{D}_{\mathcal{A}_c}\boldsymbol{F}_{\mathcal{A}_c}(t)\boldsymbol{E}_{\mathcal{A}_c} = diag\left(\boldsymbol{D}_{\mathcal{A}_{c1}}, \dots, \boldsymbol{D}_{\mathcal{A}_{cN}}\right) diag\left(\boldsymbol{F}_{\mathcal{A}_{c1}}(t), \dots, \boldsymbol{F}_{\mathcal{A}_{cN}}(t)\right) diag\left(\boldsymbol{E}_{\mathcal{A}_{c1}}, \dots, \boldsymbol{E}_{\mathcal{A}_{cN}}\right)\\
&\boldsymbol{D}_{\mathcal{B}_c}\boldsymbol{F}_{\mathcal{B}_c}(t)\boldsymbol{E}_{\mathcal{B}_c} = diag\left(\boldsymbol{D}_{\mathcal{B}_{c1}}, \dots, \boldsymbol{D}_{\mathcal{B}_{cN}}\right) diag\left(\boldsymbol{F}_{\mathcal{B}_{c1}}(t), \dots, \boldsymbol{F}_{\mathcal{B}_{cN}}(t)\right) diag\left(\boldsymbol{E}_{\mathcal{B}_{c1}}, \dots, \boldsymbol{E}_{\mathcal{B}_{cN}}\right)\\
&\boldsymbol{D}_{\mathcal{C}_c}\boldsymbol{F}_{\mathcal{C}_c}(t)\boldsymbol{E}_{\mathcal{C}_c} = diag\left(\boldsymbol{D}_{\mathcal{C}_{c1}}, \dots, \boldsymbol{D}_{\mathcal{C}_{cN}}\right) diag\left(\boldsymbol{F}_{\mathcal{C}_{c1}}(t), \dots, \boldsymbol{F}_{\mathcal{C}_{cN}}(t)\right) diag\left(\boldsymbol{E}_{\mathcal{C}_{c1}}, \dots, \boldsymbol{E}_{\mathcal{C}_{cN}}\right)\\
&\boldsymbol{D}_{\mathcal{D}_c}\boldsymbol{F}_{\mathcal{D}_c}(t)\boldsymbol{E}_{\mathcal{D}_c} = diag\left(\boldsymbol{D}_{\mathcal{D}_{c1}}, \dots, \boldsymbol{D}_{\mathcal{D}_{cN}}\right) diag\left(\boldsymbol{F}_{\mathcal{D}_{c1}}(t), \dots, \boldsymbol{F}_{\mathcal{D}_{cN}}(t)\right) diag\left(\boldsymbol{E}_{\mathcal{D}_{c1}}, \dots, \boldsymbol{E}_{\mathcal{D}_{cN}}\right)
\end{aligned}
\tag{20}
$$

uncertainty matrices in controller parameters. To simplify *(18)* we can represent the control protocol in the following form

$$
\begin{aligned}
&D^q\boldsymbol{x}_c(t) = \left(\boldsymbol{\mathcal{A}}_c + \boldsymbol{D}_{\mathcal{A}_c}\boldsymbol{F}_{\mathcal{A}_c}(t)\boldsymbol{E}_{\mathcal{A}_c}\right)\boldsymbol{x}_c(t) + \left(\boldsymbol{\mathcal{B}}_c + \boldsymbol{D}_{\mathcal{B}_c}\boldsymbol{F}_{\mathcal{B}_c}(t)\boldsymbol{E}_{\mathcal{B}_c}\right)\boldsymbol{L}_p\boldsymbol{\mathcal{Y}}(t),\\
&\boldsymbol{u} = \left(\boldsymbol{\mathcal{C}}_c + \boldsymbol{D}_{\mathcal{C}_c}\boldsymbol{F}_{\mathcal{C}_c}(t)\boldsymbol{E}_{\mathcal{C}_c}\right)\boldsymbol{x}_c + \left(\boldsymbol{\mathcal{D}}_c + \boldsymbol{D}_{\mathcal{D}_c}\boldsymbol{F}_{\mathcal{D}_c}(t)\boldsymbol{E}_{\mathcal{D}_c}\right)\boldsymbol{L}_p\boldsymbol{\mathcal{Y}}(t),\\
&\boldsymbol{x}_c(0) = \boldsymbol{x}_{c0}
\end{aligned}
\tag{21}
$$

where $\boldsymbol{L}_p = \boldsymbol{L} \otimes \boldsymbol{I}_p$, in which $\boldsymbol{L}$ is the Laplacian matrix of corresponding graph of FOMAS. Implementing control protocol *(21)* to the main system *(15)* yields

$$
\begin{aligned}
&D^q\begin{bmatrix}\boldsymbol{x}(t)\\ \boldsymbol{x}_c(t)\end{bmatrix} = \begin{bmatrix}\overline{\boldsymbol{\mathcal{A}}_{\Delta,N}} + \boldsymbol{\mathcal{B}}\left(\boldsymbol{\mathcal{D}}_c + \boldsymbol{D}_{\mathcal{D}_c}\boldsymbol{F}_{\mathcal{D}_c}(t)\boldsymbol{E}_{\mathcal{D}_c}\right)\boldsymbol{L}_p\boldsymbol{\mathcal{C}}_N & \boldsymbol{\mathcal{B}}\left(\boldsymbol{\mathcal{C}}_c + \boldsymbol{D}_{\mathcal{C}_c}\boldsymbol{F}_{\mathcal{C}_c}(t)\boldsymbol{E}_{\mathcal{C}_c}\right)\\ \left(\boldsymbol{\mathcal{B}}_c + \boldsymbol{D}_{\mathcal{B}_c}\boldsymbol{F}_{\mathcal{B}_c}(t)\boldsymbol{E}_{\mathcal{B}_c}\right)\boldsymbol{L}_p\boldsymbol{\mathcal{C}}_N & \boldsymbol{\mathcal{A}}_c + \boldsymbol{D}_{\mathcal{A}_c}\boldsymbol{F}_{\mathcal{A}_c}(t)\boldsymbol{E}_{\mathcal{A}_c}\end{bmatrix}\begin{bmatrix}\boldsymbol{x}(t)\\ \boldsymbol{x}_c(t)\end{bmatrix} +\\
&\begin{bmatrix}\boldsymbol{\phi}\left(\boldsymbol{x}(t), \boldsymbol{u}(t)\right)\\ 0\end{bmatrix} = \left(\begin{bmatrix}\boldsymbol{\mathcal{A}}_N + \boldsymbol{\mathcal{B}}\boldsymbol{\mathcal{D}}_c\boldsymbol{L}_p\boldsymbol{\mathcal{C}}_N & \boldsymbol{\mathcal{B}}\boldsymbol{\mathcal{C}}_c\\ \boldsymbol{\mathcal{B}}_c\boldsymbol{L}_p\boldsymbol{\mathcal{C}}_N & \boldsymbol{\mathcal{A}}_c\end{bmatrix} + \begin{bmatrix}\boldsymbol{\mathcal{B}}\boldsymbol{D}_{\mathcal{D}_c}\boldsymbol{F}_{\mathcal{D}_c}(t)\boldsymbol{E}_{\mathcal{D}_c}\boldsymbol{L}_p\boldsymbol{\mathcal{C}}_N & \boldsymbol{\mathcal{B}}\boldsymbol{D}_{\mathcal{C}_c}\boldsymbol{F}_{\mathcal{C}_c}(t)\boldsymbol{E}_{\mathcal{C}_c}\\ \boldsymbol{D}_{\mathcal{B}_c}\boldsymbol{F}_{\mathcal{B}_c}(t)\boldsymbol{E}_{\mathcal{B}_c}\boldsymbol{L}_p\boldsymbol{\mathcal{C}}_N & \boldsymbol{D}_{\mathcal{A}_c}\boldsymbol{F}_{\mathcal{A}_c}(t)\boldsymbol{E}_{\mathcal{A}_c}\end{bmatrix} +\right.\\
&\left.\begin{bmatrix}\boldsymbol{\Delta}\boldsymbol{\mathcal{A}}_N & 0\\ 0 & 0\end{bmatrix}\right)\begin{bmatrix}\boldsymbol{x}(t)\\ \boldsymbol{x}_c(t)\end{bmatrix} + \begin{bmatrix}\boldsymbol{\phi}\left(\boldsymbol{x}(t), \boldsymbol{u}(t)\right)\\ 0\end{bmatrix},\\
&\begin{bmatrix}\boldsymbol{x}(0)\\ \boldsymbol{x}_c(0)\end{bmatrix} = \begin{bmatrix}\boldsymbol{x}_0\\ \boldsymbol{x}_{c0}\end{bmatrix}
\end{aligned}
\tag{22}
$$

***Lemma 5.*** *[46]* $L_qC_N = C_NL_n$.

***Proof.*** $L_qC_N = \left(L \otimes I_q\right)\left(I_N \otimes \tilde{C}\right) = (LI_N) \otimes \left(I_q\tilde{C}\right) = (I_NL) \otimes \left(\tilde{C}I_n\right) = \left(I_N \otimes \tilde{C}\right)(L \otimes I_n) = C_NL_n$.

***Lemma 6.*** *[46]* $L_nA_N = A_NL_n$.

***Proof.*** $L_nA_N = (L \otimes I_n)\left(I_N \otimes \tilde{A}\right) = (LI_N) \otimes \left(I_n\tilde{A}\right) = (I_NL) \otimes \left(\tilde{A}I_n\right) = \left(I_N \otimes \tilde{A}\right)(L \otimes I_n) = A_NL_n$.

The idea is to convert the consensus problem of *(22)* into a stabilization one. In order to guarantee the consensus of FOMAS *(18)* via stabilization problem, we define a new system with the following states

$$
\overline{\boldsymbol{x}} = \boldsymbol{L}_n\boldsymbol{x}, \quad \boldsymbol{L}_n = \boldsymbol{L} \otimes \boldsymbol{I}_n. \tag{23}
$$

***Remark 5.*** *The obtained* $\overline{\boldsymbol{x}}$ *in equation (23) and its relationship with* $\boldsymbol{x}$ *is conveying the fact that it plays a relative state difference in the FOMAS. Put simply, according to the Laplacian matrix nature and equation (23), the consensus defined in* (16) *will be achieved if* $\overline{\boldsymbol{x}}$ *converges to zero. The new parameter* $\overline{\boldsymbol{x}}$ *has solved a challenge, but the second challenge has arisen by second part of the equation (23). The rank deficiency of* $\boldsymbol{L}_n$ *will result in a rank-deficient objective function* [47]. *Since most inequality constrained optimization algorithms are not able to deal with rank-deficient objective functions this is a challenge that should be overcome. A simple and effective way to overcome this obstacle is to remove* *one of the $\boldsymbol{L}$ rows and name it $\hat{\boldsymbol{L}}$. The new matrix $\hat{\boldsymbol{L}}$ has the independent rows of $\boldsymbol{L}$ and eliminates system redundancy. In the following we carry out the control scheme by means of $\hat{\boldsymbol{L}}$.*

It is obvious that $\bar{\boldsymbol{x}}$ provides a relation between $\boldsymbol{x}_i$ and its neighbors. The stability of the transformed system ensures the convergence of relative states difference of system *(18)* to zero which is equivalent to the consensus of system *(18)* defined in *(16)*. Nevertheless, the transformed system by (23) has a redundancy due to rank deficiency of $\boldsymbol{L}_n$. Removing a row of the matrix $\boldsymbol{L}$ eliminates the system redundancy, so we modify (23) as follows

$$\boldsymbol{x}_r = \hat{\boldsymbol{L}}_n \boldsymbol{x}, \quad \hat{\boldsymbol{L}}_n = \hat{\boldsymbol{L}} \otimes \boldsymbol{I}_n. \tag{24}$$

In order to express the pseudo state space representation of system with respect to $\boldsymbol{x}_r$, consider

$$\begin{aligned} &\boldsymbol{\mathcal{A}}_{\Delta,N} = \boldsymbol{I}_N \otimes \boldsymbol{\Delta}\widetilde{\boldsymbol{\mathcal{A}}}, \\ &\boldsymbol{\Delta}\widetilde{\boldsymbol{\mathcal{A}}} = \widetilde{\boldsymbol{\mathcal{R}}}\delta(\sigma)\widetilde{\boldsymbol{\mathcal{N}}} \end{aligned} \tag{25}$$

where

$$\begin{aligned} &\delta = \{\delta_j \mid \|\delta_j\| > \|\delta_i\|, i \neq j\}, \;\; i = 1, \dots, N, \\ &\delta(\sigma) = \boldsymbol{\mathcal{Z}}(\sigma)[\boldsymbol{I} + \boldsymbol{\mathcal{J}}\boldsymbol{\mathcal{Z}}(\sigma)]^{-1} \,. \end{aligned} \tag{26}$$

Multiplying both sides of Equation (22) by $diag\left(\hat{\boldsymbol{L}}_n, \boldsymbol{I}_{N.n_c}\right)$ and then, after some calculations we obtain the closed-loop system

$$\begin{aligned} &D^q \begin{bmatrix} \boldsymbol{x}_r(t) \\ \boldsymbol{x}_c(t) \end{bmatrix} = \\ &\begin{bmatrix} \hat{\boldsymbol{L}}_n \boldsymbol{\mathcal{A}}_{\Delta,N} \hat{\boldsymbol{L}}_n^{\uparrow} + \hat{\boldsymbol{L}}_n \boldsymbol{\mathcal{B}}\left(\boldsymbol{\mathcal{D}}_c + \boldsymbol{D}_{\boldsymbol{\mathcal{D}}_c} \boldsymbol{F}_{\boldsymbol{\mathcal{D}}_c}(t) \boldsymbol{E}_{\boldsymbol{\mathcal{D}}_c}\right) \boldsymbol{L}_p \boldsymbol{\mathcal{C}}_N \hat{\boldsymbol{L}}_n^{\uparrow} & \hat{\boldsymbol{L}}_n \boldsymbol{\mathcal{B}}\left(\boldsymbol{\mathcal{C}}_c + \boldsymbol{D}_{\boldsymbol{\mathcal{C}}_c} \boldsymbol{F}_{\boldsymbol{\mathcal{C}}_c}(t) \boldsymbol{E}_{\boldsymbol{\mathcal{C}}_c}\right) \\ \left(\boldsymbol{\mathcal{B}}_c + \boldsymbol{D}_{\boldsymbol{\mathcal{B}}_c} \boldsymbol{F}_{\boldsymbol{\mathcal{B}}_c}(t) \boldsymbol{E}_{\boldsymbol{\mathcal{B}}_c}\right) \boldsymbol{L}_p \boldsymbol{\mathcal{C}}_N \hat{\boldsymbol{L}}_n^{\uparrow} & \boldsymbol{\mathcal{A}}_c + \boldsymbol{D}_{\boldsymbol{\mathcal{A}}_c} \boldsymbol{F}_{\boldsymbol{\mathcal{A}}_c}(t) \boldsymbol{E}_{\boldsymbol{\mathcal{A}}_c} \end{bmatrix} \begin{bmatrix} \boldsymbol{x}_r(t) \\ \boldsymbol{x}_c(t) \end{bmatrix} + \\ &\begin{bmatrix} \hat{\boldsymbol{L}}_n \boldsymbol{\phi}\left(\boldsymbol{x}(t), \boldsymbol{u}(t)\right) \\ 0 \end{bmatrix}. \end{aligned} \tag{27}$$

Given that $\boldsymbol{L}_q \boldsymbol{C}_N = \boldsymbol{C}_N \boldsymbol{L}_n, \boldsymbol{L}_n \boldsymbol{A}_N = \boldsymbol{A}_N \boldsymbol{L}_n$ [46], to match the matrices dimensions in the resulting system matrix, it can be easily obtained that $\hat{\boldsymbol{L}}_n \boldsymbol{\mathcal{A}}_{\Delta,N} \hat{\boldsymbol{L}}_n^{\uparrow} = \boldsymbol{I}_{N-1} \otimes \left(\widetilde{\boldsymbol{\mathcal{A}}} + \boldsymbol{\Delta}\widetilde{\boldsymbol{\mathcal{A}}}\right) = \boldsymbol{\mathcal{A}}_{\Delta,N-1}$ and $\boldsymbol{L}_p \boldsymbol{\mathcal{C}}_N \hat{\boldsymbol{L}}_n^{\uparrow} = \boldsymbol{\mathcal{C}}_N \boldsymbol{L}_n \boldsymbol{L}_n^{\uparrow}$. Then, defining $\boldsymbol{\mathcal{C}}_r = \boldsymbol{\mathcal{C}}_N \boldsymbol{L}_n \hat{\boldsymbol{L}}_n^{\uparrow}$ and $\boldsymbol{\phi}_{\boldsymbol{r}}\left(\boldsymbol{x}_{\mathbf{r}}(t), \boldsymbol{u}(t)\right) = \hat{\boldsymbol{L}}_n \boldsymbol{\phi}\left(\boldsymbol{x}(t), \boldsymbol{u}(t)\right)$, and $\boldsymbol{\mathcal{X}}(t) = [\boldsymbol{x}_r^T(t) \quad \boldsymbol{x}_c^T(t)]^T$ the closed-loop system is achieved as follows

$$\begin{aligned} &D^q \boldsymbol{\mathcal{X}}(t) = \boldsymbol{\Phi}(\boldsymbol{\mathcal{X}}(t), t) = \boldsymbol{\mathcal{A}}_{cl,\Delta} \boldsymbol{\mathcal{X}}(t) + \begin{bmatrix} \boldsymbol{\phi}_r\left(\boldsymbol{x}_r(t), \boldsymbol{u}(t)\right) \\ 0 \end{bmatrix}, \\ &\boldsymbol{\mathcal{X}}(0) = \boldsymbol{\mathcal{X}}_0 = [\boldsymbol{x}_{r0}^T \quad \boldsymbol{x}_{c0}^T]^T \end{aligned} \tag{28}$$

where

$$\begin{aligned} &\boldsymbol{\mathcal{A}}_{cl,\Delta} = \boldsymbol{\mathcal{A}}_{\psi} + \boldsymbol{\mathcal{A}}_{\Delta}, \quad \boldsymbol{\mathcal{X}}(t) = [\boldsymbol{x}_r^T(t) \quad \boldsymbol{x}_c^T(t)]^T, \\ &\boldsymbol{\mathcal{A}}_{\psi} = \begin{bmatrix} \boldsymbol{\mathcal{A}}_{N-1} + \hat{\boldsymbol{L}}_n \boldsymbol{\mathcal{B}} \boldsymbol{\mathcal{D}}_c \boldsymbol{\mathcal{C}}_r & \hat{\boldsymbol{L}}_n \boldsymbol{\mathcal{B}} \boldsymbol{\mathcal{C}}_c \\ \boldsymbol{\mathcal{B}}_c \boldsymbol{\mathcal{C}}_r & \boldsymbol{\mathcal{A}}_c \end{bmatrix}, \\ &\boldsymbol{\mathcal{A}}_{\Delta} = \begin{bmatrix} \boldsymbol{\Delta}\boldsymbol{\mathcal{A}}_{N-1} + \hat{\boldsymbol{L}}_n \boldsymbol{\mathcal{B}} \boldsymbol{D}_{\boldsymbol{\mathcal{D}}_c} \boldsymbol{F}_{\boldsymbol{\mathcal{D}}_c}(t) \boldsymbol{E}_{\boldsymbol{\mathcal{D}}_c} \boldsymbol{\mathcal{C}}_r & \hat{\boldsymbol{L}}_n \boldsymbol{\mathcal{B}} \boldsymbol{D}_{\boldsymbol{\mathcal{C}}_c} \boldsymbol{F}_{\boldsymbol{\mathcal{C}}_c}(t) \boldsymbol{E}_{\boldsymbol{\mathcal{C}}_c} \\ \boldsymbol{D}_{\boldsymbol{\mathcal{B}}_c} \boldsymbol{F}_{\boldsymbol{\mathcal{B}}_c}(t) \boldsymbol{E}_{\boldsymbol{\mathcal{B}}_c} \boldsymbol{\mathcal{C}}_r & \boldsymbol{D}_{\boldsymbol{\mathcal{A}}_c} \boldsymbol{F}_{\boldsymbol{\mathcal{A}}_c}(t) \boldsymbol{E}_{\boldsymbol{\mathcal{A}}_c} \end{bmatrix}, \end{aligned} \tag{29}$$

and

$$\begin{aligned} &\boldsymbol{\Delta}\boldsymbol{\mathcal{A}}_{N-1} = \boldsymbol{\mathcal{R}} \boldsymbol{\Delta}(\sigma) \boldsymbol{\mathcal{N}} \\ &\boldsymbol{\mathcal{R}} = \boldsymbol{I}_{N-1} \otimes \widetilde{\boldsymbol{\mathcal{R}}}, \;\; \boldsymbol{\mathcal{N}} = \boldsymbol{I}_{N-1} \otimes \widetilde{\boldsymbol{\mathcal{N}}}, \\ &\widehat{\boldsymbol{\mathcal{Z}}}(\sigma) = \boldsymbol{I}_{N-1} \otimes \boldsymbol{\mathcal{Z}}(\sigma), \;\; \widehat{\boldsymbol{\mathcal{J}}} = I_{N-1} \otimes \boldsymbol{\mathcal{J}}, \\ &\boldsymbol{\Delta}(\sigma) = \boldsymbol{I}_{N-1} \otimes \delta(\sigma) = \widehat{\boldsymbol{\mathcal{Z}}}(\sigma)\left[\boldsymbol{I} - \widehat{\boldsymbol{\mathcal{J}}} \widehat{\boldsymbol{\mathcal{Z}}}(\sigma)\right]^{-1} \end{aligned} \tag{30}$$

**Theorem 1** The output feedback controller *(18)* solves the consensus problem of multi-agent system *(1)* If there exist positive constants $\tau_i$ for $i = 1, \dots, 5, \mu, \xi$ and positive definite matrix $\boldsymbol{P}$ such that the following matrix inequality holds:

$$\begin{bmatrix} \boldsymbol{P}\boldsymbol{\mathcal{A}}_{\psi} + \boldsymbol{\mathcal{A}}_{\psi}^{T}\boldsymbol{P} + \boldsymbol{\Pi} & \boldsymbol{P}\begin{bmatrix} \hat{\boldsymbol{L}}_n \boldsymbol{\mathcal{B}} \boldsymbol{D}_{\mathcal{D}_c} \\ \boldsymbol{0} \end{bmatrix} & \boldsymbol{P}\begin{bmatrix} \hat{\boldsymbol{L}}_n \boldsymbol{\mathcal{B}} \boldsymbol{D}_{\mathcal{C}_c} \\ \boldsymbol{0} \end{bmatrix} & \boldsymbol{P}\begin{bmatrix} \boldsymbol{I} \\ \boldsymbol{0} \end{bmatrix} & \boldsymbol{P}\begin{bmatrix} \boldsymbol{0} \\ \boldsymbol{D}_{\mathcal{B}_c} \end{bmatrix} & \boldsymbol{P}\begin{bmatrix} \boldsymbol{0} \\ \boldsymbol{D}_{\mathcal{A}_c} \end{bmatrix} & \begin{bmatrix} \boldsymbol{\mathcal{R}} & \boldsymbol{P}^{T}\boldsymbol{\mathcal{N}}^{T} \end{bmatrix} \\ \star & -\tau_1 \boldsymbol{I} & \boldsymbol{0} & \boldsymbol{0} & \boldsymbol{0} & \boldsymbol{0} & \boldsymbol{0} \\ \star & \star & -\tau_2 \boldsymbol{I} & \boldsymbol{0} & \boldsymbol{0} & \boldsymbol{0} & \boldsymbol{0} \\ \star & \star & \star & -\tau_3 \boldsymbol{I} & \boldsymbol{0} & \boldsymbol{0} & \boldsymbol{0} \\ \star & \star & \star & \star & -\tau_4 \boldsymbol{I} & \boldsymbol{0} & \boldsymbol{0} \\ \star & \star & \star & \star & \star & -\tau_5 \boldsymbol{I} & \boldsymbol{0} \\ \star & \star & \star & \star & \star & \star & \begin{bmatrix} \mu \boldsymbol{I} & -\mu \boldsymbol{I} \\ -\mu \boldsymbol{I} & \widehat{\boldsymbol{\mathcal{W}}} + \mu \boldsymbol{I} \end{bmatrix} \end{bmatrix} < 0 \qquad (31)$$

where

$$\boldsymbol{\Pi} = \begin{bmatrix} \left(\left\|\boldsymbol{E}_{\mathcal{D}_c}\boldsymbol{\mathcal{C}}_r\right\| + \left\|\boldsymbol{E}_{\mathcal{B}_c}\boldsymbol{\mathcal{C}}_r\right\|\right)\boldsymbol{I} + \xi \boldsymbol{I} & \boldsymbol{0} \\ \star & \left(\left\|\boldsymbol{E}_{\mathcal{C}_c}\right\| + \left\|\boldsymbol{E}_{\mathcal{A}_c}\right\|\right)\boldsymbol{I} + \xi \boldsymbol{I} \end{bmatrix}, \widehat{\boldsymbol{\mathcal{W}}} = Sym(\hat{\boldsymbol{\mathcal{J}}}), \text{ in which } \hat{\boldsymbol{\mathcal{J}}} \text{ is provided in (30).}$$

***Proof.*** The stability of transformed system (28) guarantees the consensus of (1). Considering the closed-loop system (28), for any $\boldsymbol{\mathcal{X}}_1(t) = [\boldsymbol{x}_{r1}^{T}(t) \quad \boldsymbol{x}_{c1}^{T}(t)]^{T}$ and $\boldsymbol{\mathcal{X}}_2(t) = [\boldsymbol{x}_{r2}^{T}(t) \quad \boldsymbol{x}_{c2}^{T}(t)]^{T}$ we have

$$\left\| \boldsymbol{\mathcal{A}}_{cl,\Delta}\boldsymbol{\mathcal{X}}_1(t) + \begin{bmatrix} \boldsymbol{\phi}_r\left(\boldsymbol{x}_{r_1}(t), \boldsymbol{u}(t)\right) \\ 0 \end{bmatrix} - \boldsymbol{\mathcal{A}}_{cl,\Delta}\boldsymbol{\mathcal{X}}_2(t) - \begin{bmatrix} \boldsymbol{\phi}_r\left(\boldsymbol{x}_{r2}(t), \boldsymbol{u}(t)\right) \\ 0 \end{bmatrix} \right\|_2 \leq \left\|\boldsymbol{\mathcal{A}}_{cl,\Delta}\right\|_2 \|\boldsymbol{\mathcal{X}}_1(t) - \boldsymbol{\mathcal{X}}_2(t)\|_2 + \left\|\hat{\boldsymbol{L}}_n\right\|_2 \left\|\boldsymbol{\phi}\left(\boldsymbol{x}_{r1}(t), \boldsymbol{u}(t)\right) - \boldsymbol{\phi}\left(\boldsymbol{x}_{r2}(t), \boldsymbol{u}(t)\right)\right\|_2 \qquad (32)$$

then According to $\phi\left(x_r(t), u(t)\right) = \left[\tilde{\phi}^{T}\left(x_{r1}(t), u_1(t)\right), \dots, \tilde{\phi}^{T}\left(x_{rN}(t), u_N(t)\right)\right]^{T}$ one can obtain

$$\left\|\hat{\boldsymbol{L}}_n\right\|_2 \left\|\boldsymbol{\phi}\left(\boldsymbol{x}_{r1}(t), \boldsymbol{u}(t)\right) - \boldsymbol{\phi}\left(\boldsymbol{x}_{r2}(t), \boldsymbol{u}(t)\right)\right\|_2 = \left\|\hat{\boldsymbol{L}}_n\right\|_2 \left\| \begin{matrix} \tilde{\boldsymbol{\phi}}\left(\boldsymbol{x}_{r11}(t), \boldsymbol{u}(t)\right) - \tilde{\boldsymbol{\phi}}\left(\boldsymbol{x}_{r12}(t), \boldsymbol{u}(t)\right) \\ \vdots \\ \tilde{\boldsymbol{\phi}}\left(\boldsymbol{x}_{rN1}(t), \boldsymbol{u}(t)\right) - \tilde{\boldsymbol{\phi}}\left(\boldsymbol{x}_{rN2}(t), \boldsymbol{u}(t)\right) \end{matrix} \right\|_2 \qquad (33)$$

since $\tilde{\boldsymbol{\phi}}(x, u)$ satisfies the Assumption 3, Lipschitz condition implies that

$$\left\|\hat{\boldsymbol{L}}_n\right\|_2 \left\| \begin{matrix} \tilde{\boldsymbol{\phi}}\left(\boldsymbol{x}_{r11}(t), \boldsymbol{u}(t)\right) - \tilde{\boldsymbol{\phi}}\left(\boldsymbol{x}_{r12}(t), \boldsymbol{u}(t)\right) \\ \vdots \\ \tilde{\boldsymbol{\phi}}\left(\boldsymbol{x}_{rN1}(t), \boldsymbol{u}(t)\right) - \tilde{\boldsymbol{\phi}}\left(\boldsymbol{x}_{rN2}(t), \boldsymbol{u}(t)\right) \end{matrix} \right\|_2 < \xi_1 \left\|\hat{\boldsymbol{L}}_n\right\|_2 \left\| \begin{bmatrix} \boldsymbol{x}_{r11}(t) - \boldsymbol{x}_{r12}(t) \\ \vdots \\ \boldsymbol{x}_{rN1}(t) - \boldsymbol{x}_{rN2}(t) \end{bmatrix} \right\|_2 \qquad (34)$$

since for every matrix like $Q \in \mathbb{R}^{N \times L}$ according to the Frobenius norm definition $\|Q\|_2 = \sqrt{\sum_{i=1}^{N} \sum_{j=1}^{L} q_{ij}^2}$ for the (34) it can be easily obtained

$$\xi_1 \left\|\hat{\boldsymbol{L}}_n\right\|_2 \left\| \begin{bmatrix} \boldsymbol{x}_{r11}(t) - \boldsymbol{x}_{r12}(t) \\ \vdots \\ \boldsymbol{x}_{rN1}(t) - \boldsymbol{x}_{rN2}(t) \end{bmatrix} \right\|_2 \leq \xi_1 \left\|\hat{\boldsymbol{L}}_n\right\|_2 \left\| \begin{bmatrix} \boldsymbol{x}_{r11}(t) - \boldsymbol{x}_{r12}(t) \\ \vdots \\ \boldsymbol{x}_{rN1}(t) - \boldsymbol{x}_{rN2}(t) \\ \boldsymbol{x}_{c11}(t) - \boldsymbol{x}_{c12}(t) \\ \vdots \\ \boldsymbol{x}_{cN1}(t) - \boldsymbol{x}_{cN2}(t) \end{bmatrix} \right\|_2 = \xi_1 \left\|\hat{\boldsymbol{L}}_n\right\|_2 \|\boldsymbol{\mathcal{X}}_1(t) - \boldsymbol{\mathcal{X}}_2(t)\|_2 \qquad (35)$$

Matrices $\boldsymbol{\mathcal{A}}_{cl,\Delta}$ and $\hat{\boldsymbol{L}}_n$ have bounded elements there exist constants $M_1 > 0$ and $M_2 > 0$ such that $\left\|\hat{\boldsymbol{L}}_n\right\|_2 \leq M_1$ and $\left\|\hat{\boldsymbol{L}}_n\right\|_2 \leq M_2$. Substituting *(35)* into inequality *(32)*, it implies that

$$\left\| \boldsymbol{\mathcal{A}}_{cl,\Delta}\boldsymbol{\mathcal{X}}_1(t) + \begin{bmatrix} \boldsymbol{\phi}_r\left(\boldsymbol{x}_{r_1}(t), \boldsymbol{u}(t)\right) \\ 0 \end{bmatrix} - \boldsymbol{\mathcal{A}}_{cl,\Delta}\boldsymbol{\mathcal{X}}_2(t) - \begin{bmatrix} \boldsymbol{\phi}_r\left(\boldsymbol{x}_{r2}(t), \boldsymbol{u}(t)\right) \\ 0 \end{bmatrix} \right\|_2 \leq (M_1 + \xi_1 M_2)\|\boldsymbol{\mathcal{X}}_1(t) - \boldsymbol{\mathcal{X}}_2(t)\|_2, \qquad (36)$$

this yields that $\boldsymbol{\Phi}(\boldsymbol{\mathcal{X}}(t), t)$ is Lipschitz in $\boldsymbol{\mathcal{X}}(t)$.

Define $\boldsymbol{\Phi}_{\boldsymbol{\mathcal{X}}}(\boldsymbol{\mathcal{X}}(t),t) = \boldsymbol{\mathcal{A}}_{cl,\Delta}\boldsymbol{\mathcal{X}}(t) + \begin{bmatrix}\boldsymbol{\phi}_r\big(\boldsymbol{x}_r(t),\boldsymbol{u}(t)\big)\\ 0\end{bmatrix}$ a continuous function mapping from a set $\mathbb{R}_\epsilon = \{(t,\boldsymbol{\mathcal{X}}): 0 \le t \le a \ and \ \|\boldsymbol{\mathcal{X}}-\boldsymbol{\mathcal{X}}_{\boldsymbol{0}}\| \le b\}$ to $\mathbb{R}^{N(n+n_c)}$. $\boldsymbol{\Phi}(\boldsymbol{\mathcal{X}}(t),t)$ is bounded on $\mathbb{R}_\epsilon$ with upper bound $M_3 > 0$. It follows from Lemma 2 and Lemma 3 that, the solution of *(28)* is given by

$$\boldsymbol{\mathcal{X}}(t) = \boldsymbol{\mathcal{X}}_*\big(t^q/\Gamma(q+1)\big) \qquad (37)$$

where $\boldsymbol{\mathcal{X}}_*(v)$ is the solution of following integer order differential equation

$$\begin{aligned}&\boldsymbol{d\mathcal{X}}_*(v)/\boldsymbol{d}v = \boldsymbol{\mathcal{A}}_{cl,\Delta}\boldsymbol{\mathcal{X}}_*(v) + \Xi\big(\boldsymbol{\mathcal{X}}_*(v),\boldsymbol{u}_*(v)\big),\\ &\boldsymbol{\mathcal{X}}_*(0) = [\boldsymbol{x}_{r0}^T \quad \boldsymbol{x}_{c0}^T]^T\end{aligned} \qquad (38)$$

with

$$\begin{aligned}&\boldsymbol{\mathcal{X}}_*(v) = \boldsymbol{\mathcal{X}}\left(t - \left(t^q - v\Gamma(q+1)\right)^{1/q}\right),\\ &\boldsymbol{x}_{r*}(v) = \boldsymbol{x}_r\left(t - \left(t^q - v\Gamma(q+1)\right)^{1/q}\right),\\ &\boldsymbol{x}_{c*}(v) = \boldsymbol{x}_c\left(t - \left(t^q - v\Gamma(q+1)\right)^{1/q}\right),\\ &\boldsymbol{u}_*(v) = \boldsymbol{u}\left(t - \left(t^q - v\Gamma(q+1)\right)^{1/q}\right),\\ &\Xi\big(\boldsymbol{\mathcal{X}}_*(t),\boldsymbol{u}_*(t)\big) = \begin{bmatrix}\boldsymbol{\phi}_r\big(\boldsymbol{x}_{r*}(t),\boldsymbol{u}_*(t)\big)\\ 0\end{bmatrix}\end{aligned} \qquad (39)$$

If the system *(38)* be stable, it guarantees the system *(28)* stability and consequently, the consensus of the system *(1)*. Consider a candidate Lyapunov function for *(38)* as follows

$$\boldsymbol{V}(v) = \boldsymbol{\mathcal{X}}_*^T(v)\boldsymbol{P}\boldsymbol{\mathcal{X}}_*(v) \qquad (40)$$

where $\boldsymbol{P}$ is a symmetric positive definite matrix. Then its time derivative is calculated as

$$\begin{aligned}&\boldsymbol{dV}(v)/\boldsymbol{d}v = \dot{\boldsymbol{\mathcal{X}}}_*^T(v)\boldsymbol{P}\boldsymbol{\mathcal{X}}_*(v) + \boldsymbol{\mathcal{X}}_*^T(v)\boldsymbol{P}\dot{\boldsymbol{\mathcal{X}}}_*(v) = \left(\boldsymbol{\mathcal{A}}_{cl,\Delta}\boldsymbol{\mathcal{X}}_*(v) + \Xi\big(\boldsymbol{\mathcal{X}}_*(v),\boldsymbol{u}_*(v)\big)\right)^T \boldsymbol{P}\boldsymbol{\mathcal{X}}_*(v) +\\ &\boldsymbol{\mathcal{X}}_*^T(v)\boldsymbol{P}\left(\boldsymbol{\mathcal{A}}_{cl,\Delta}\boldsymbol{\mathcal{X}}_*(v) + \Xi\big(\boldsymbol{\mathcal{X}}_*(v),\boldsymbol{u}_*(v)\big)\right) = \boldsymbol{\mathcal{X}}_*^T(v)\big(\boldsymbol{P}\boldsymbol{\mathcal{A}}_\psi + \boldsymbol{\mathcal{A}}_\psi^T\boldsymbol{P} + sym\{\boldsymbol{P}\boldsymbol{\mathcal{R}}\boldsymbol{\Delta}(\sigma)\boldsymbol{\mathcal{N}}\}\big)\boldsymbol{\mathcal{X}}_*(v) +\\ &\boldsymbol{\mathcal{X}}_*^T(v)\left(\boldsymbol{P}\begin{bmatrix}\hat{\boldsymbol{L}}_n\boldsymbol{\mathcal{B}}\boldsymbol{D}_{\boldsymbol{\mathcal{D}}_c}\boldsymbol{F}_{\boldsymbol{\mathcal{D}}_c}(t)\boldsymbol{E}_{\boldsymbol{\mathcal{D}}_c}\boldsymbol{\mathcal{C}}_r & \hat{\boldsymbol{L}}_n\boldsymbol{\mathcal{B}}\boldsymbol{D}_{\boldsymbol{\mathcal{C}}_c}\boldsymbol{F}_{\boldsymbol{\mathcal{C}}_c}(t)\boldsymbol{E}_{\boldsymbol{\mathcal{C}}_c}\\ \boldsymbol{D}_{\boldsymbol{\mathcal{B}}_c}\boldsymbol{F}_{\boldsymbol{\mathcal{B}}_c}(t)\boldsymbol{E}_{\boldsymbol{\mathcal{B}}_c}\boldsymbol{\mathcal{C}}_r & \boldsymbol{D}_{\boldsymbol{\mathcal{A}}_c}\boldsymbol{F}_{\boldsymbol{\mathcal{A}}_c}(t)\boldsymbol{E}_{\boldsymbol{\mathcal{A}}_c}\end{bmatrix}\boldsymbol{\mathcal{X}}_*(v) + \Xi\big(\boldsymbol{\mathcal{X}}_*(v),\boldsymbol{u}_*(v)\big)\right) +\\ &\left(\boldsymbol{P}\begin{bmatrix}\hat{\boldsymbol{L}}_n\boldsymbol{\mathcal{B}}\boldsymbol{D}_{\boldsymbol{\mathcal{D}}_c}\boldsymbol{F}_{\boldsymbol{\mathcal{D}}_c}(t)\boldsymbol{E}_{\boldsymbol{\mathcal{D}}_c}\boldsymbol{\mathcal{C}}_r & \hat{\boldsymbol{L}}_n\boldsymbol{\mathcal{B}}\boldsymbol{D}_{\boldsymbol{\mathcal{C}}_c}\boldsymbol{F}_{\boldsymbol{\mathcal{C}}_c}(t)\boldsymbol{E}_{\boldsymbol{\mathcal{C}}_c}\\ \boldsymbol{D}_{\boldsymbol{\mathcal{B}}_c}\boldsymbol{F}_{\boldsymbol{\mathcal{B}}_c}(t)\boldsymbol{E}_{\boldsymbol{\mathcal{B}}_c}\boldsymbol{\mathcal{C}}_r & \boldsymbol{D}_{\boldsymbol{\mathcal{A}}_c}\boldsymbol{F}_{\boldsymbol{\mathcal{A}}_c}(t)\boldsymbol{E}_{\boldsymbol{\mathcal{A}}_c}\end{bmatrix}\boldsymbol{\mathcal{X}}_*(v) + \Xi\big(\boldsymbol{\mathcal{X}}_*(v),\boldsymbol{u}_*(v)\big)\right)^T \boldsymbol{\mathcal{X}}_*(v)\end{aligned} \qquad (41)$$

the equation *(41)* can be rewritten as

$$\begin{aligned}&\boldsymbol{dV}(v)/\boldsymbol{d}v = \boldsymbol{\mathcal{X}}_*^T(v)\big(\boldsymbol{P}\boldsymbol{\mathcal{A}}_\psi + \boldsymbol{\mathcal{A}}_\psi^T\boldsymbol{P} + sym\{\boldsymbol{P}\boldsymbol{\mathcal{R}}\boldsymbol{\Delta}(\sigma)\boldsymbol{\mathcal{N}}\}\big)\boldsymbol{\mathcal{X}}_*(v) +\\ &\boldsymbol{\mathcal{X}}_*^T(v)\boldsymbol{P}\left(\begin{bmatrix}\hat{\boldsymbol{L}}_n\boldsymbol{\mathcal{B}}\boldsymbol{D}_{\boldsymbol{\mathcal{D}}_c}z_1 + \hat{\boldsymbol{L}}_n\boldsymbol{\mathcal{B}}\boldsymbol{D}_{\boldsymbol{\mathcal{C}}_c}z_2 + z_3\\ \boldsymbol{D}_{\boldsymbol{\mathcal{B}}_c}z_4 + \boldsymbol{D}_{\boldsymbol{\mathcal{A}}_c}z_5\end{bmatrix}\right) + \left(\begin{bmatrix}\hat{\boldsymbol{L}}_n\boldsymbol{\mathcal{B}}\boldsymbol{D}_{\boldsymbol{\mathcal{D}}_c}z_1 + \hat{\boldsymbol{L}}_n\boldsymbol{\mathcal{B}}\boldsymbol{D}_{\boldsymbol{\mathcal{C}}_c}z_2 + z_3\\ \boldsymbol{D}_{\boldsymbol{\mathcal{B}}_c}z_4 + \boldsymbol{D}_{\boldsymbol{\mathcal{A}}_c}z_5\end{bmatrix}\right)^T \boldsymbol{P}\boldsymbol{\mathcal{X}}_*(v)\end{aligned} \qquad (42)$$

where $z_1 = \boldsymbol{F}_{\boldsymbol{\mathcal{D}}_c}(t)\boldsymbol{E}_{\boldsymbol{\mathcal{D}}_c}\boldsymbol{\mathcal{C}}_r\boldsymbol{x}_{r*}$, $z_2 = \boldsymbol{F}_{\boldsymbol{\mathcal{C}}_c}(t)\boldsymbol{E}_{\boldsymbol{\mathcal{C}}_c}\boldsymbol{x}_{c*}(v)$, $z_3 = \boldsymbol{\phi}_r\big(\boldsymbol{x}_{r*}(t),\boldsymbol{u}_*(t)\big)$, $z_4 = \boldsymbol{F}_{\boldsymbol{\mathcal{B}}_c}(t)\boldsymbol{E}_{\boldsymbol{\mathcal{B}}_c}\boldsymbol{\mathcal{C}}_r\boldsymbol{x}_{r*}(v)$, $z_5 = \boldsymbol{F}_{\boldsymbol{\mathcal{A}}_c}(t)\boldsymbol{E}_{\boldsymbol{\mathcal{A}}_c}\boldsymbol{x}_{c*}(v)$. Introducing vector $Z = [\boldsymbol{\mathcal{X}}_*^T(v)\ z_1^T\ z_2^T\ z_3^T\ z_4^T\ z_5^T]^T$, the Equation (42) rearranged

$$\frac{\boldsymbol{dV}(v)}{\boldsymbol{d}v} = Z^T\begin{bmatrix}\boldsymbol{P}\boldsymbol{\mathcal{A}}_\psi + \boldsymbol{\mathcal{A}}_\psi^T\boldsymbol{P} + sym\{\boldsymbol{P}\boldsymbol{\mathcal{R}}\boldsymbol{\Delta}(\sigma)\boldsymbol{\mathcal{N}}\} & \boldsymbol{P}\begin{bmatrix}\hat{\boldsymbol{L}}_n\boldsymbol{\mathcal{B}}\boldsymbol{D}_{\boldsymbol{\mathcal{D}}_c}\\ \boldsymbol{0}\end{bmatrix} & \boldsymbol{P}\begin{bmatrix}\hat{\boldsymbol{L}}_n\boldsymbol{\mathcal{B}}\boldsymbol{D}_{\boldsymbol{\mathcal{C}}_c}\\ \boldsymbol{0}\end{bmatrix} & \boldsymbol{P}\begin{bmatrix}\boldsymbol{I}\\ \boldsymbol{0}\end{bmatrix} & \boldsymbol{P}\begin{bmatrix}\boldsymbol{0}\\ \boldsymbol{D}_{\boldsymbol{\mathcal{B}}_c}\end{bmatrix} & \boldsymbol{P}\begin{bmatrix}\boldsymbol{0}\\ \boldsymbol{D}_{\boldsymbol{\mathcal{A}}_c}\end{bmatrix}\\ \star & \boldsymbol{0} & \boldsymbol{0} & \boldsymbol{0} & \boldsymbol{0} & \boldsymbol{0}\\ \star & \star & \boldsymbol{0} & \boldsymbol{0} & \boldsymbol{0} & \boldsymbol{0}\\ \star & \star & \star & \boldsymbol{0} & \boldsymbol{0} & \boldsymbol{0}\\ \star & \star & \star & \star & \boldsymbol{0} & \boldsymbol{0}\\ \star & \star & \star & \star & \star & \boldsymbol{0}\end{bmatrix}Z \qquad (43)$$

According to the direct Lyapunov approach, the stability conditions for the system *(38)* are $\boldsymbol{V}(\nu) > 0$ and $\boldsymbol{dV}(\nu)/\boldsymbol{d}\nu < 0$. Equation *(40)* shows that $\boldsymbol{V}(\nu)$ is positive, and the second condition holds if $\boldsymbol{dV}(\nu)/\boldsymbol{d}\nu$, defined in (43) be negative.

***Remark 6.*** *In continuous time systems without a delay, using quadratic functions as Lyapunov candidate is common. Desirable results of using Quadratic Functions as a Lyapunov candidate can be found in literature* [14], [39], [48]. *According to the literature we selected the quadratic function (40) and got good and non-conservative results which is discussed in the simulation section.*

It follows from *(12)* and *(19)*, that

$$\begin{aligned} & z_1^T z_1 < \left\|\boldsymbol{E}_{\mathcal{D}_c}\boldsymbol{\mathcal{C}}_r\right\| \boldsymbol{x}_{r*}^T(\nu)\boldsymbol{x}_{r*}(\nu), z_2^T z_2 < \left\|\boldsymbol{E}_{\mathcal{C}_c}\right\| \boldsymbol{x}_{c*}^T(\nu)\boldsymbol{x}_{c*}(\nu), z_3^T z_3 < \xi_1^2 \boldsymbol{\mathcal{X}}_*^T(\nu)\boldsymbol{\mathcal{X}}_*(\nu), \\ & z_4^T z_4 < \left\|\boldsymbol{E}_{\mathcal{B}_c}\boldsymbol{\mathcal{C}}_r\right\| \boldsymbol{x}_{r*}^T(\nu)\boldsymbol{x}_{r*}(\nu), z_5^T z_5 < \left\|\boldsymbol{E}_{\mathcal{A}_c}\right\| \boldsymbol{x}_{c*}^T(\nu)\boldsymbol{x}_{c*}(\nu) \end{aligned} \tag{44}$$

The combination of inequalities *(44)* with respect to $\boldsymbol{Z}$, yields

$$Z^T \begin{bmatrix} \begin{bmatrix} \left(\left\|\boldsymbol{E}_{\mathcal{D}_c}\boldsymbol{\mathcal{C}}_r\right\| + \left\|\boldsymbol{E}_{\mathcal{B}_c}\boldsymbol{\mathcal{C}}_r\right\|\right)\boldsymbol{I} + \xi_1^2\boldsymbol{I} & \boldsymbol{0} \\ \star & \left(\left\|\boldsymbol{E}_{\mathcal{C}_c}\right\| + \left\|\boldsymbol{E}_{\mathcal{A}_c}\right\|\right)\boldsymbol{I} + \xi_1^2\boldsymbol{I} \end{bmatrix} & \boldsymbol{0} & \boldsymbol{0} & \boldsymbol{0} & \boldsymbol{0} & \boldsymbol{0} \\ \star & -\boldsymbol{I} & \boldsymbol{0} & \boldsymbol{0} & \boldsymbol{0} & \boldsymbol{0} \\ \star & \star & -\boldsymbol{I} & \boldsymbol{0} & \boldsymbol{0} & \boldsymbol{0} \\ \star & \star & \star & -\boldsymbol{I} & \boldsymbol{0} & \boldsymbol{0} \\ \star & \star & \star & \star & -\boldsymbol{I} & \boldsymbol{0} \\ \star & \star & \star & \star & \star & -\boldsymbol{I} \end{bmatrix} Z < 0 \tag{45}$$

Applying S-Procedure on $\boldsymbol{dV}(\nu)/\boldsymbol{d}\nu < 0$, where $\boldsymbol{dV}(\nu)/\boldsymbol{d}\nu$ is defined in (43), and *(45)*, and also defining $\xi = \xi_1^2$, it can be obtained that

$$Z^T \boldsymbol{\Sigma} Z < 0 \tag{46}$$

where

$$\boldsymbol{\Sigma} = \begin{bmatrix} \boldsymbol{P}\boldsymbol{\mathcal{A}}_\psi + \boldsymbol{\mathcal{A}}_\psi^T\boldsymbol{P} + sym\{\boldsymbol{P}\boldsymbol{\mathcal{R}}\boldsymbol{\Delta}(\sigma)\boldsymbol{\mathcal{N}}\} + \boldsymbol{\Pi} & \boldsymbol{P}\begin{bmatrix} \hat{\boldsymbol{L}}_n\boldsymbol{\mathcal{B}}\boldsymbol{D}_{\mathcal{D}_c} \\ \boldsymbol{0} \end{bmatrix} & \boldsymbol{P}\begin{bmatrix} \hat{\boldsymbol{L}}_n\boldsymbol{\mathcal{B}}\boldsymbol{D}_{\mathcal{C}_c} \\ \boldsymbol{0} \end{bmatrix} & \boldsymbol{P}\begin{bmatrix} \boldsymbol{I} \\ \boldsymbol{0} \end{bmatrix} & \boldsymbol{P}\begin{bmatrix} \boldsymbol{0} \\ \boldsymbol{D}_{\mathcal{B}_c} \end{bmatrix} & \boldsymbol{P}\begin{bmatrix} \boldsymbol{0} \\ \boldsymbol{D}_{\mathcal{A}_c} \end{bmatrix} \\ \star & -\boldsymbol{\tau}_1\boldsymbol{I} & \boldsymbol{0} & \boldsymbol{0} & \boldsymbol{0} & \boldsymbol{0} \\ \star & \star & -\boldsymbol{\tau}_2\boldsymbol{I} & \boldsymbol{0} & \boldsymbol{0} & \boldsymbol{0} \\ \star & \star & \star & -\boldsymbol{\tau}_3\boldsymbol{I} & \boldsymbol{0} & \boldsymbol{0} \\ \star & \star & \star & \star & -\boldsymbol{\tau}_4\boldsymbol{I} & \boldsymbol{0} \\ \star & \star & \star & \star & \star & -\boldsymbol{\tau}_5\boldsymbol{I} \end{bmatrix} \tag{47}$$

Define

$$\begin{aligned} & \widehat{\boldsymbol{\mathcal{W}}} = Sym\left(\widehat{\boldsymbol{\mathcal{J}}}\right), \\ & \widehat{\boldsymbol{\mathcal{Q}}} = \widehat{\boldsymbol{\mathcal{W}}}^{-1/2}\left(\boldsymbol{\mathcal{N}} + \boldsymbol{\mathcal{R}}^T\boldsymbol{P}^T\right) - \widehat{\boldsymbol{\mathcal{W}}}^{-1/2}\boldsymbol{\Delta}(\sigma)\boldsymbol{\mathcal{N}}. \end{aligned} \tag{48}$$

it follows from Lemma 4 that $Sym\{\boldsymbol{\Delta}(\sigma)\} - \boldsymbol{\Delta}^T(\sigma)\widehat{\boldsymbol{\mathcal{W}}}\boldsymbol{\Delta}(\sigma) > 0$, and the following inequality holds

$$\begin{aligned} & -\widehat{\boldsymbol{\mathcal{Q}}}^T\widehat{\boldsymbol{\mathcal{Q}}} \leq \boldsymbol{0} \Leftrightarrow -\boldsymbol{sym}\left\{\boldsymbol{\mathcal{N}}^T\widehat{\boldsymbol{\mathcal{W}}}^{-1}\boldsymbol{\mathcal{R}}^T\boldsymbol{P}^T\right\} - \boldsymbol{P}\boldsymbol{\mathcal{R}}\boldsymbol{\Delta}(\sigma)\boldsymbol{\mathcal{R}}^T\boldsymbol{P}^T - \boldsymbol{\mathcal{N}}^T\widehat{\boldsymbol{\mathcal{W}}}^{-1}\boldsymbol{\mathcal{N}} + \boldsymbol{\mathcal{N}}^T\left(sym\{\boldsymbol{\Delta}(\sigma)\} - \right. \\ & \left. \boldsymbol{\Delta}^T(\sigma)\widehat{\boldsymbol{\mathcal{W}}}\boldsymbol{\Delta}(\sigma)\right)\boldsymbol{\mathcal{N}} + sym\{\boldsymbol{P}\boldsymbol{\mathcal{R}}\boldsymbol{\Delta}(\sigma)\boldsymbol{\mathcal{N}}\} \leq 0\,. \end{aligned} \tag{49}$$

It implies from (49) that

$$sym\{\boldsymbol{P}\boldsymbol{\mathcal{R}}\boldsymbol{\Delta}(\sigma)\boldsymbol{\mathcal{N}}\} \leq \boldsymbol{sym}\left\{\boldsymbol{\mathcal{N}}^T\widehat{\boldsymbol{\mathcal{W}}}^{-1}\boldsymbol{\mathcal{R}}^T\boldsymbol{P}^T\right\} + \boldsymbol{P}\boldsymbol{\mathcal{R}}\boldsymbol{\Delta}(\sigma)\boldsymbol{\mathcal{R}}^T\boldsymbol{P}^T + \boldsymbol{\mathcal{N}}^T\widehat{\boldsymbol{\mathcal{W}}}^{-1}\boldsymbol{\mathcal{N}}\,. \tag{50}$$

Inequality *(50)* is equivalent to that there exist and $\mu > 0$ such that

$$\boldsymbol{sym}\{\boldsymbol{P}\boldsymbol{\mathcal{R}}\boldsymbol{\Delta}(\boldsymbol{\sigma})\boldsymbol{\mathcal{N}}\} \leq \boldsymbol{sym}\left\{\boldsymbol{\mathcal{N}}^T\widehat{\boldsymbol{\mathcal{W}}}^{-1}\boldsymbol{\mathcal{R}}^T\boldsymbol{P}^T\right\} + \boldsymbol{P}\boldsymbol{\mathcal{R}}\boldsymbol{\Delta}(\sigma)\boldsymbol{\mathcal{R}}^T\boldsymbol{P}^T + \boldsymbol{\mathcal{N}}^T\left(\widehat{\boldsymbol{\mathcal{W}}}^{-1} + \mu^{-1}\boldsymbol{I}\right)\boldsymbol{\mathcal{N}}\,, \tag{51}$$

which is equivalent to that there exist and $\mu > 0$ such that

$$
\boldsymbol{sym}\{\boldsymbol{P}\boldsymbol{\mathcal{R}}\boldsymbol{\Delta}(\boldsymbol{\sigma})\boldsymbol{\mathcal{N}}\} \le [\boldsymbol{\mathcal{R}} \quad \boldsymbol{P}^T\boldsymbol{\mathcal{N}}^T] \begin{bmatrix} \widehat{\boldsymbol{\mathcal{W}}}^{-1} + \mu^{-1} I & \widehat{\boldsymbol{\mathcal{W}}}^{-1} \\ \widehat{\boldsymbol{\mathcal{W}}}^{-1} & \widehat{\boldsymbol{\mathcal{W}}}^{-1} \end{bmatrix} \begin{bmatrix} \boldsymbol{\mathcal{R}}^T \\ \boldsymbol{\mathcal{N}}\boldsymbol{P} \end{bmatrix} =
$$

$$
[\boldsymbol{\mathcal{R}} \quad \boldsymbol{P}^T\boldsymbol{\mathcal{N}}^T] \begin{bmatrix} \mu \boldsymbol{I} & -\mu \boldsymbol{I} \\ -\mu \boldsymbol{I} & \widehat{\boldsymbol{\mathcal{W}}} + \mu \boldsymbol{I} \end{bmatrix}^{-1} \begin{bmatrix} \boldsymbol{\mathcal{R}}^T \\ \boldsymbol{\mathcal{N}}\boldsymbol{P} \end{bmatrix}. \tag{52}
$$

Substituting *(52)* into inequality *(47)*, and applying Schur complement completes the proof. ■

***Remark 7.*** *As observed in Theorem 1, the inequality (31)* $\boldsymbol{\mathcal{A}}_{\psi}$ *contains some varying parameters multiplied by* $\boldsymbol{P}$ *where the inequity (31) can be described in the general BMI form* $F(x,y) = F_0 + \sum_{i=1}^{m} x_i F_i + \sum_{j=1}^{n} y_i G_i + \sum_{i=1}^{m} \sum_{j=1}^{n} x_i y_i H_{ij} > 0$ *where* $F_i, G_j, H_{ij} \in \mathbb{R}^{n\times n}$ *are all fixed symmetric matrices and* $x \in \mathbb{R}^m, y \in \mathbb{R}^n$ *are variables* [49]. *A BMI is convex in* $\boldsymbol{x}$ *for fixed* $\boldsymbol{y}$ *and convex in* $\boldsymbol{y}$ *for fixed* $\boldsymbol{x}$. *Nonconvex sets can be easily described by BMIs. As a result, much wider phenomenon's constraints can be described and modeled by BMIs, such as wide range of optimization and control problems. But the main problem dealing with BMIs is its computation difficulties than LMIs* [50], [51]. *On the other hand, an LMI has the form* $\boldsymbol{F}(\boldsymbol{x}) = \boldsymbol{F}_0 + \sum_{i=1}^{m} \boldsymbol{x}_i \boldsymbol{F}_i > 0$ *where* $\boldsymbol{x} \in \mathbb{R}^m$, $\boldsymbol{F}_i \in \mathbb{R}^{n\times n}$ [49]. *Given that an LMI is a convex constraint, it can be solved easier than a BMI with different solvers. That is why LMIs are preferred in control theory to design controllers. As a result, in the following theorem, we are going to describe the resulted inequality in Theorem 1, in the form of an LMI.*

**Theorem 2** The output feedback controller *(18)* solves the consensus problem of the system *(1)* with $0 < \alpha < 1$, if there exist positive constants $\tau_i$ $(i = 1, \dots, 5)$, $\mu$, $\xi$ and positive definite matrice1s $\overline{\boldsymbol{p}}_u \in \mathbb{R}^{n\times n}$, $\boldsymbol{p}_{di} \in \mathbb{R}^{\boldsymbol{n_c}\times \boldsymbol{n_c}}(i = 1, \dots, N)$ and matrices $\mathfrak{A} = diag(\mathfrak{a}_1, \dots, \mathfrak{a}_N), \mathfrak{B} = diag(\mathfrak{b}_1, \dots, \mathfrak{b}_N)$ and $\mathfrak{c}_i$, $\mathfrak{d}_i$ for $i = 1, \dots, N$, such that the following matrix inequality holds:

$$
\begin{bmatrix}
\boldsymbol{\Pi} + \boldsymbol{\Omega} & \boldsymbol{P}\begin{bmatrix} \hat{\boldsymbol{L}}_n \boldsymbol{\mathcal{B}} \boldsymbol{D}_{\mathcal{D}_c} \\ \boldsymbol{0} \end{bmatrix} & \boldsymbol{P}\begin{bmatrix} \hat{\boldsymbol{L}}_n \boldsymbol{\mathcal{B}} \boldsymbol{D}_{\mathcal{C}_c} \\ \boldsymbol{0} \end{bmatrix} & \boldsymbol{P}\begin{bmatrix} \boldsymbol{I} \\ \boldsymbol{0} \end{bmatrix} & \boldsymbol{P}\begin{bmatrix} \boldsymbol{0} \\ \boldsymbol{D}_{\mathcal{B}_c} \end{bmatrix} & \boldsymbol{P}\begin{bmatrix} \boldsymbol{0} \\ \boldsymbol{D}_{\mathcal{A}_c} \end{bmatrix} & [\boldsymbol{\mathcal{R}} \quad \boldsymbol{P}\boldsymbol{\mathcal{N}}^T] \\
\star & -\boldsymbol{\tau}_1 \boldsymbol{I} & \boldsymbol{0} & \boldsymbol{0} & \boldsymbol{0} & \boldsymbol{0} & \boldsymbol{0} \\
\star & \star & -\boldsymbol{\tau}_2 \boldsymbol{I} & \boldsymbol{0} & \boldsymbol{0} & \boldsymbol{0} & \boldsymbol{0} \\
\star & \star & \star & -\boldsymbol{\tau}_3 \boldsymbol{I} & \boldsymbol{0} & \boldsymbol{0} & \boldsymbol{0} \\
\star & \star & \star & \star & -\boldsymbol{\tau}_4 \boldsymbol{I} & \boldsymbol{0} & \boldsymbol{0} \\
\star & \star & \star & \star & \star & -\boldsymbol{\tau}_5 \boldsymbol{I} & \boldsymbol{0} \\
\star & \star & \star & \star & \star & \star & \begin{bmatrix} \mu \boldsymbol{I} & -\mu \boldsymbol{I} \\ -\mu \boldsymbol{I} & \widehat{\boldsymbol{\mathcal{W}}} + \mu \boldsymbol{I} \end{bmatrix}
\end{bmatrix} < 0 \tag{53}
$$

where

$$
\boldsymbol{\Pi} = \begin{bmatrix} \left(\|\boldsymbol{E}_{\mathcal{D}_c} \boldsymbol{\mathcal{C}}_r\| + \|\boldsymbol{E}_{\mathcal{B}_c} \boldsymbol{\mathcal{C}}_r\|\right) \boldsymbol{I} + \xi \boldsymbol{I} & \boldsymbol{0} \\ \star & \left(\|\boldsymbol{E}_{\mathcal{C}_c}\| + \|\boldsymbol{E}_{\mathcal{A}_c}\|\right) \boldsymbol{I} + \xi \boldsymbol{I} \end{bmatrix},
$$

$$
\boldsymbol{\Omega} = \begin{bmatrix} \boldsymbol{P}_u \boldsymbol{\mathcal{A}}_{N-1} + \boldsymbol{\mathcal{A}}_{N-1}^T \boldsymbol{P}_u + \mathfrak{D} \boldsymbol{\mathcal{C}}_r + \boldsymbol{\mathcal{C}}_r^T \mathfrak{D}^T & \mathfrak{C} + \boldsymbol{\mathcal{C}}_r^T \mathfrak{B}^T \\ \star & \mathfrak{A} + \mathfrak{A}^T \end{bmatrix},
$$

$$
\mathfrak{D} = \left(\hat{\boldsymbol{L}} \otimes \boldsymbol{1}_{\boldsymbol{n}}\right) \circ \left(\boldsymbol{1}_{\boldsymbol{N-1}} \otimes [\mathfrak{d}_1, \dots, \mathfrak{d}_N]\right), \mathfrak{C} = \left(\hat{\boldsymbol{L}} \otimes \boldsymbol{J}_{\boldsymbol{n}\times \boldsymbol{nc}}\right) \circ \left(\boldsymbol{1}_{\boldsymbol{N-1}} \otimes [\mathfrak{c}_1, \dots, \mathfrak{c}_N]\right)
$$

$$
\boldsymbol{P}_u = \boldsymbol{I}_{N-1} \otimes \overline{\boldsymbol{p}}_u,\ \boldsymbol{P}_d = \boldsymbol{diag}(\boldsymbol{p}_{d1}, \dots, \boldsymbol{p}_{dN}),\ \boldsymbol{P} = \begin{bmatrix} \boldsymbol{P}_u & \boldsymbol{0} \\ \boldsymbol{0} & \boldsymbol{P}_d \end{bmatrix}
$$

$\widehat{\boldsymbol{\mathcal{W}}} = Sym(\hat{\mathfrak{J}})$, in which $\hat{\mathfrak{J}}$ is provided in *(30)*. Moreover, the controller matrices $\boldsymbol{\mathcal{A}}_c$, $\boldsymbol{\mathcal{B}}_c$, $\boldsymbol{\mathcal{C}}_c$ and $\boldsymbol{\mathcal{D}}_c$ are as follows

$$
\begin{gathered}
\boldsymbol{\mathcal{A}}_c = \boldsymbol{P}_d^{-1} \mathfrak{A},\quad \boldsymbol{\mathcal{B}}_c = \boldsymbol{P}_d^{-1} \mathfrak{B}, \\
\boldsymbol{\mathcal{C}}_c = diag\left(\boldsymbol{B}_1^{\uparrow} \overline{\boldsymbol{p}}_u^{-1} \mathfrak{c}_1, \dots, \boldsymbol{B}_N^{\uparrow} \overline{\boldsymbol{p}}_u^{-1} \mathfrak{c}_N\right), \boldsymbol{\mathcal{D}}_c = diag\left(\boldsymbol{B}_1^{\uparrow} \overline{\boldsymbol{p}}_u^{-1} \mathfrak{d}_1, \dots, \boldsymbol{B}_N^{\uparrow} \overline{\boldsymbol{p}}_u^{-1} \mathfrak{d}_N\right)
\end{gathered} \tag{54}
$$

***Proof.*** According to the proof of Theorem 1 , the output feedback controller *(18)* solves the consensus problem of the system *(1)* if the inequality *(31)* holds. To deal with multiplication of variables, according to $\boldsymbol{P}$ defined by

$$\boldsymbol{P}=\begin{bmatrix}\boldsymbol{P}_u & \mathbf{0}\\ \mathbf{0} & \boldsymbol{P}_d\end{bmatrix}, \boldsymbol{P}_u=\boldsymbol{I}_{N-1}\otimes\overline{\boldsymbol{p}}_u, \boldsymbol{P}_d=diag\left(\boldsymbol{P}_{d_1},\dots,\boldsymbol{P}_{d_N}\right). \tag{55}$$

we expand the matrix $\boldsymbol{P}\boldsymbol{\mathcal{A}}_\psi+\boldsymbol{\mathcal{A}}_\psi^T\boldsymbol{P}$

$$\boldsymbol{P}\boldsymbol{\mathcal{A}}_\psi+\boldsymbol{\mathcal{A}}_\psi^T\boldsymbol{P}=\begin{bmatrix}\boldsymbol{P}_u\boldsymbol{\mathcal{A}}_{N-1}+\boldsymbol{\mathcal{A}}_{N-1}\boldsymbol{P}_u+\boldsymbol{P}_u\hat{\boldsymbol{L}}_n\boldsymbol{\mathcal{B}}\boldsymbol{\mathcal{D}}_c\boldsymbol{\mathcal{C}}_r+\boldsymbol{\mathcal{C}}_r^T\boldsymbol{\mathcal{D}}_c^T\boldsymbol{\mathcal{B}}^T\hat{\boldsymbol{L}}_n^T\boldsymbol{P}_u & \boldsymbol{P}_u\hat{\boldsymbol{L}}_n\boldsymbol{\mathcal{B}}\boldsymbol{\mathcal{C}}_c+\boldsymbol{\mathcal{C}}_r^T\boldsymbol{\mathcal{B}}_c^T\boldsymbol{P}_d\\ \boldsymbol{P}_d\boldsymbol{\mathcal{B}}_c\boldsymbol{\mathcal{C}}_r+\boldsymbol{\mathcal{C}}_c^T\boldsymbol{\mathcal{B}}^T\hat{\boldsymbol{L}}_n^T\boldsymbol{P}_u & \boldsymbol{P}_d\boldsymbol{\mathcal{A}}_c+\boldsymbol{\mathcal{A}}_c^T\boldsymbol{P}_d\end{bmatrix}$$

$$\boldsymbol{P}_u\hat{\boldsymbol{L}}_n\boldsymbol{\mathcal{B}}\boldsymbol{\mathcal{D}}_c=\left(\hat{\boldsymbol{L}}\otimes\boldsymbol{J}_{\boldsymbol{n}\times\boldsymbol{p}}\right)\circ\left(\mathbf{1}_{N-1}\otimes[\overline{\boldsymbol{p}}_u\boldsymbol{B}_1\overline{\boldsymbol{\mathcal{D}}}_{c1},\dots,\overline{\boldsymbol{p}}_u\boldsymbol{B}_N\overline{\boldsymbol{\mathcal{D}}}_{cN}]\right)$$

$$\boldsymbol{P}_u\hat{\boldsymbol{L}}_n\boldsymbol{\mathcal{B}}\boldsymbol{\mathcal{C}}_c=\left(\hat{\boldsymbol{L}}\otimes\boldsymbol{J}_{\boldsymbol{n}\times\boldsymbol{n_c}}\right)\circ\left(\mathbf{1}_{\boldsymbol{N-1}}\otimes[\overline{\boldsymbol{p}}_u\boldsymbol{B}_1\overline{\boldsymbol{\mathcal{C}}}_{c1},\dots,\overline{\boldsymbol{p}}_u\boldsymbol{B}_N\overline{\boldsymbol{\mathcal{C}}}_{cN}]\right).$$

Now, the following change of the variables completes the proof

$$\mathfrak{A}=\boldsymbol{P}_d\boldsymbol{\mathcal{A}}_c,\ \ \mathfrak{B}=\boldsymbol{P}_d\boldsymbol{\mathcal{B}}_c,\ \ \mathfrak{c}_i=\overline{\boldsymbol{p}}_u\boldsymbol{B}_i\overline{\boldsymbol{\mathcal{C}}}_{ci},\ \ \mathfrak{d}_i=\overline{\boldsymbol{p}}_u\boldsymbol{B}_i\overline{\boldsymbol{\mathcal{D}}}_{ci}.\ \blacksquare$$

***Corollary 1*** *Consider the nonlinear fractional-order multi-agent system (1) without nonlinear term, the output dynamic controller solves the consensus problem if there exist two real symmetric positive definite matrices* $\overline{\boldsymbol{p}}_u\in\mathbb{R}^{n\times n}$, $\boldsymbol{p}_{di}\in\mathbb{R}^{\boldsymbol{n_c}\times\boldsymbol{n_c}}(i=1,\dots,N)$, *and matrices* $\mathfrak{A}=diag(\mathfrak{a}_1,\dots,\mathfrak{a}_N),\mathfrak{C}=diag(\mathfrak{c}_1,\dots,\mathfrak{c}_N)$ *and* $\mathfrak{b}_i$, $\mathfrak{d}_i$ *for* $i=1,\dots,N$ *and a real constant* $\mu>0$ *such that*

$$\begin{bmatrix}\sin\theta\left(\boldsymbol{\mathcal{A}}_{\mathfrak{C}}\boldsymbol{P}+\boldsymbol{P}\boldsymbol{\mathcal{A}}_{\mathfrak{C}}^T\right) & \sin\theta\,\widehat{\boldsymbol{\mathcal{R}}} & I_2\otimes\left(\boldsymbol{P}^T\widehat{\boldsymbol{\mathcal{N}}}^T\right)\\ \star & -\mu\boldsymbol{I} & \mu\boldsymbol{I}\\ \star & \star & -\widehat{\boldsymbol{W}}-\mu\boldsymbol{I}\end{bmatrix}<0, \tag{56}$$

*where* $\widehat{\boldsymbol{W}}=Sym(\widehat{\boldsymbol{\mathfrak{J}}})$, *in which* $\widehat{\boldsymbol{\mathfrak{J}}}$ *is provided in (30), and*

$$\boldsymbol{P}=\begin{bmatrix}\boldsymbol{P}_u & \mathbf{0}\\ \mathbf{0} & \boldsymbol{p}_d\end{bmatrix}=\begin{bmatrix}\boldsymbol{I}_{N-1}\otimes\overline{\boldsymbol{p}}_u & \mathbf{0}_{(N-1).n,N.n_c}\\ \mathbf{0}_{N.n_c,(N-1).n} & diag(\boldsymbol{p}_{d1},\dots,\boldsymbol{p}_{dN})\end{bmatrix},$$

$$\boldsymbol{\mathcal{A}}_{\mathfrak{C}}=\begin{bmatrix}\boldsymbol{\mathcal{A}}_{N-1}\boldsymbol{P}_u+\hat{\boldsymbol{L}}_n\boldsymbol{\mathcal{B}}\mathfrak{D} & \hat{\boldsymbol{L}}_n\boldsymbol{\mathcal{B}}\mathfrak{C}\\ \mathfrak{B} & \mathfrak{A}\end{bmatrix},\ \overline{\boldsymbol{\mathcal{R}}}=[\boldsymbol{\mathcal{R}}^T\quad\mathbf{0}]^T,\ \overline{\boldsymbol{\mathcal{N}}}=[\boldsymbol{\mathcal{N}}\quad\mathbf{0}], \tag{57}$$

$$\mathfrak{B}=\begin{bmatrix}diag(\mathfrak{b}_1,\dots,\mathfrak{b}_{N-1})\\ -\mathbf{1}_{N-1}^T\otimes\mathfrak{b}_N\end{bmatrix},\mathfrak{D}=\begin{bmatrix}diag(\mathfrak{d}_1,\dots,\mathfrak{d}_{N-1})\\ -\mathbf{1}_{N-1}^T\otimes\mathfrak{d}_N\end{bmatrix}$$

$$\theta=q\pi/2.$$

*The controller matrices* $\boldsymbol{\mathcal{A}}_c$, $\boldsymbol{\mathcal{B}}_c$, $\boldsymbol{\mathcal{C}}_c$ *and* $\boldsymbol{\mathcal{D}}_c$ *can be obtain as follows*

$$\begin{aligned}&\boldsymbol{\mathcal{A}}_c=\mathfrak{A}{\boldsymbol{P}_d}^{-1},\boldsymbol{\mathcal{B}}_c=diag\left(\mathfrak{b}_1\overline{\boldsymbol{p}}_u^{-1}\widetilde{\boldsymbol{\mathcal{C}}}^{\dagger},\dots,\mathfrak{b}_N\overline{\boldsymbol{p}}_u^{-1}\widetilde{\boldsymbol{\mathcal{C}}}^{\dagger}\right)\\ &\boldsymbol{\mathcal{C}}_c=\mathfrak{C}{\boldsymbol{P}_d}^{-1},\boldsymbol{\mathcal{D}}_c=diag\left(\mathfrak{d}_1\overline{\boldsymbol{p}}_u^{-1}\widetilde{\boldsymbol{\mathcal{C}}}^{\dagger},\dots,\mathfrak{d}_N\overline{\boldsymbol{p}}_u^{-1}\widetilde{\boldsymbol{\mathcal{C}}}^{\dagger}\right)\end{aligned} \tag{58}$$

***Proof.*** By eliminating the nonlinear term, the closed-loop system *(28)* becomes

$$D^q\boldsymbol{\mathcal{X}}(t)=\boldsymbol{\mathcal{A}}_{cl,\Delta}\boldsymbol{\mathcal{X}}(t) \tag{59}$$

where

$$\begin{aligned}&\boldsymbol{\mathcal{A}}_{cl,\Delta}=\boldsymbol{\mathcal{A}}_\psi+\boldsymbol{\mathcal{A}}_\Delta,\\ &\boldsymbol{\mathcal{A}}_\psi=\begin{bmatrix}\boldsymbol{\mathcal{A}}_{N-1}+\hat{\boldsymbol{L}}_n\boldsymbol{\mathcal{B}}\boldsymbol{\mathcal{D}}_c\boldsymbol{\mathcal{C}}_r & \hat{\boldsymbol{L}}_n\boldsymbol{\mathcal{B}}\boldsymbol{\mathcal{C}}_c\\ \boldsymbol{\mathcal{B}}_c\boldsymbol{\mathcal{C}}_r & \boldsymbol{\mathcal{A}}_c\end{bmatrix},\ \boldsymbol{\mathcal{A}}_\Delta=\begin{bmatrix}\boldsymbol{\mathcal{R}}\boldsymbol{\Delta}(\sigma)\boldsymbol{\mathcal{N}} & \mathbf{0}\\ \mathbf{0} & \mathbf{0}\end{bmatrix},\end{aligned} \tag{60}$$

$\boldsymbol{\mathcal{R}}$, $\boldsymbol{\Delta}(\sigma)$, and $\boldsymbol{\mathcal{N}}$ are as defined in *(30)*. According to Lemma 1, system *(59)* is asymptotically stable if and only if there exist two real symmetric positive definite matrices $\boldsymbol{P}_{k1}\in\mathbb{R}^{(N-1)n+N.n_c\times(N-1)n+N.n_c}$, $k=1,2$, and two skew-symmetric matrices $\boldsymbol{P}_{k2}\in\mathbb{R}^{(N-1)n+N.n_c\times(N-1)n+N.n_c}$, $k=1,2$, such that inequality (3) holds for $\boldsymbol{\mathcal{A}}_{cl,\Delta}$ defined in (60).

Assuming $P_{12}=P_{22}=\mathbf{0}$, $P_{11}=P_{21}=\boldsymbol{P}=diag(\boldsymbol{P}_u,\boldsymbol{P}_d)$ we can obtain $\sum_{i=1}^{2}Sym\left\{\Theta_{i1}\otimes\left(\boldsymbol{\mathcal{A}}_{cl,\Delta}\boldsymbol{P}_{i1}\right)\right\}<0\Leftrightarrow Sym\left\{\sin\theta\,\boldsymbol{\mathcal{A}}_{cl,\Delta}\boldsymbol{P}\right\}<0\Leftrightarrow Sym\left\{\sin\theta\left(\boldsymbol{\mathcal{A}}_\psi+\boldsymbol{\mathcal{A}}_\Delta\right)\boldsymbol{P}\right\}<0\Leftrightarrow\Psi+\sin\theta\,Sym\left\{\widehat{\boldsymbol{\mathcal{R}}}\boldsymbol{\Delta}(\sigma)\widehat{\boldsymbol{\mathcal{N}}}\widehat{\boldsymbol{P}}\right\}<0$ where $\Psi=Sym\left\{\sin\theta\,\boldsymbol{\mathcal{A}}_\psi\boldsymbol{P}\right\}$, $\widehat{\boldsymbol{\mathcal{R}}}=[\boldsymbol{\mathcal{R}}^T\quad\mathbf{0}]^T$, $\widehat{\boldsymbol{\mathcal{N}}}=[\boldsymbol{\mathcal{N}}\quad\mathbf{0}]$. Define $Q=\widehat{W}^{-1/2}\left(\widehat{\boldsymbol{\mathcal{R}}}^T+\widehat{\boldsymbol{\mathcal{N}}}\boldsymbol{P}\right)-\widehat{W}^{\frac{1}{2}}\boldsymbol{\Delta}^T(\sigma)\widehat{\boldsymbol{\mathcal{R}}}^T$. Similar to the proof of Theorem 1 , above discussion is equivalent to that there exist $\mu>0$ such that

$$\Psi + [\hat{\boldsymbol{\mathcal{R}}} \quad \boldsymbol{P}^T\widehat{\boldsymbol{\mathcal{N}}}^T]\begin{bmatrix} \widehat{W}^{-1}+\mu^{-1}I & \widehat{W}^{-1} \\ \widehat{W}^{-1} & \widehat{W}^{-1} \end{bmatrix}\begin{bmatrix} \hat{\boldsymbol{\mathcal{R}}}^T \\ \widehat{\boldsymbol{\mathcal{N}}}\boldsymbol{P} \end{bmatrix} = \Psi + [\hat{\boldsymbol{\mathcal{R}}} \quad \boldsymbol{P}^T\widehat{\boldsymbol{\mathcal{N}}}^T]\begin{bmatrix} \mu I & -\mu I \\ -\mu I & \widehat{W}+\mu I \end{bmatrix}^{-1}\begin{bmatrix} \hat{\boldsymbol{\mathcal{R}}}^T \\ \widehat{\boldsymbol{\mathcal{N}}}\boldsymbol{P} \end{bmatrix} < 0. \tag{61}$$

Applying Schur complement, inequality *(61)* becomes equivalent to

$$\begin{bmatrix} Sym\{\sin\theta\, \boldsymbol{\mathcal{A}}_{\psi}\boldsymbol{P}\} & \hat{\boldsymbol{\mathcal{R}}} & \boldsymbol{P}^T\widehat{\boldsymbol{\mathcal{N}}}^T \\ \star & -\mu I & \mu I \\ \star & \star & -(\widehat{W}+\mu I) \end{bmatrix} < 0\,. \tag{62}$$

Expanding $\boldsymbol{\mathcal{A}}_{\psi}\boldsymbol{P}$ leads to

$$\boldsymbol{\mathcal{A}}_{\psi}\boldsymbol{P} = \begin{bmatrix} \boldsymbol{\mathcal{A}}_{N-1}+\hat{\boldsymbol{L}}_n\boldsymbol{\mathcal{B}}\boldsymbol{\mathcal{D}}_c\boldsymbol{\mathcal{C}}_r & \hat{\boldsymbol{L}}_n\boldsymbol{\mathcal{B}}\boldsymbol{\mathcal{C}}_c \\ \boldsymbol{\mathcal{B}}_c\boldsymbol{\mathcal{C}}_r & \boldsymbol{\mathcal{A}}_c \end{bmatrix}\begin{bmatrix} \boldsymbol{P}_u & \boldsymbol{0} \\ \boldsymbol{0} & \boldsymbol{P}_d \end{bmatrix} = \begin{bmatrix} \boldsymbol{\mathcal{A}}_{N-1}\boldsymbol{P}_u+\hat{\boldsymbol{L}}_n\boldsymbol{\mathcal{B}}\boldsymbol{\mathcal{D}}_c\boldsymbol{\mathcal{C}}_r\boldsymbol{P}_u & \hat{\boldsymbol{L}}_n\boldsymbol{\mathcal{B}}\boldsymbol{\mathcal{C}}_c\boldsymbol{P}_d \\ \boldsymbol{\mathcal{B}}_c\boldsymbol{\mathcal{C}}_r\boldsymbol{P}_u & \boldsymbol{\mathcal{A}}_c\boldsymbol{P}_d \end{bmatrix}. \tag{63}$$

According to the definition of $\boldsymbol{\mathcal{C}}_r$ we get $\boldsymbol{\mathcal{D}}_c\boldsymbol{\mathcal{C}}_r\boldsymbol{P}_u = diag(\overline{\boldsymbol{\mathcal{D}}}_{c1},\dots,\overline{\boldsymbol{\mathcal{D}}}_{cN})(\boldsymbol{I}_N\otimes\widetilde{\boldsymbol{\mathcal{C}}})(\boldsymbol{L}\otimes\boldsymbol{I}_n)(\hat{\boldsymbol{L}}\otimes\boldsymbol{I}_n)^{\uparrow}(\boldsymbol{I}_{N-1}\otimes\overline{\boldsymbol{p}}_u)$, which can be rewritten as $\boldsymbol{\mathcal{D}}_c\boldsymbol{\mathcal{C}}_r\boldsymbol{P}_u = (diag(\overline{\boldsymbol{\mathcal{D}}}_{c1},\dots,\overline{\boldsymbol{\mathcal{D}}}_{cN})\boldsymbol{L}\hat{\boldsymbol{L}}^{\uparrow}\otimes\widetilde{\boldsymbol{\mathcal{C}}}\overline{\boldsymbol{p}}_u)$. It can be obtain that $\boldsymbol{L}\hat{\boldsymbol{L}}^{\uparrow} = \begin{bmatrix} \boldsymbol{I}_{N-1} \\ -\boldsymbol{1}_{N-1}^T \end{bmatrix}$ then we have $\boldsymbol{\mathcal{D}}_c\boldsymbol{\mathcal{C}}_r\boldsymbol{P}_u = \begin{bmatrix} diag(\overline{\boldsymbol{\mathcal{D}}}_{c1}\widetilde{\boldsymbol{\mathcal{C}}}\overline{\boldsymbol{p}}_u,\dots,\overline{\boldsymbol{\mathcal{D}}}_{cN-1}\widetilde{\boldsymbol{\mathcal{C}}}\overline{\boldsymbol{p}}_u) \\ -\boldsymbol{1}_{N-1}^T\otimes\overline{\boldsymbol{\mathcal{D}}}_{cN}\widetilde{\boldsymbol{\mathcal{C}}}\overline{\boldsymbol{p}}_u \end{bmatrix}$ and, with the same procedure we have $\boldsymbol{\mathcal{B}}_c\boldsymbol{\mathcal{C}}_r\boldsymbol{P}_u = \begin{bmatrix} diag(\overline{\boldsymbol{\mathcal{B}}}_{c1}\widetilde{\boldsymbol{\mathcal{C}}}\overline{\boldsymbol{p}}_u,\dots,\overline{\boldsymbol{\mathcal{B}}}_{cN-1}\widetilde{\boldsymbol{\mathcal{C}}}\overline{\boldsymbol{p}}_u) \\ -\boldsymbol{1}_{N-1}^T\otimes\overline{\boldsymbol{\mathcal{B}}}_{cN}\widetilde{\boldsymbol{\mathcal{C}}}\overline{\boldsymbol{p}}_u \end{bmatrix}$. Finally, the following change of the variables completes the proof

$$\mathfrak{A} = \boldsymbol{\mathcal{A}}_c\boldsymbol{P}_d,\ \mathfrak{b}_i = \overline{\boldsymbol{\mathcal{B}}}_{ci}\widetilde{\boldsymbol{\mathcal{C}}}\overline{\boldsymbol{p}}_u,\ \mathfrak{C} = \boldsymbol{\mathcal{C}}_c\boldsymbol{P}_d,\ \mathfrak{d}_i = \overline{\boldsymbol{\mathcal{D}}}_{ci}\widetilde{\boldsymbol{\mathcal{C}}}\overline{\boldsymbol{p}}_u,\ \blacksquare \tag{64}$$

***Corollary 2*** *In Theorem 1 and Theorem 2 the obtained dynamic output feedback controllers can be reduced to static ones by putting* $n_c = 0$. *For a system with open loop system* $A(s)$ *and the feedback system* $F(s)$, *the closed loop transfer function with an output feedback controller would be* $G(s) = A(s)/(1 + A(s)F(s))$. *The difference between static and dynamic output feedback controller is in* $F(s)$. *For a static output feedback controller* $F(s)$ *is a scalar gain but in a dynamic one it has dynamics which brings more parameters to the controller. The advantages of a dynamic output feedback controller were discussed in Remark 4.*

## IV. SIMULATION

In this section, the proposed examples demonstrate the effectiveness of the designed decentralized dynamic output feedback controllers for the consensus of fractional-order permanent magnet synchronous multi-motor velocity and, a numerical example. Various solvers and parsers can be utilized to determine variables satisfying feasibility problem. In this paper simulation results are obtained using YALMIP parser [52], implemented as a toolbox in MATLAB [53].

### A. Permanent magnet synchronous motor (PMSM)

Consider a number of motors operating together in manufacturing industries, such as textile and paper mills, which can be modeled in the form of multi-agent systems. Using multi-motor setup instead of traditional mechanical coupling, for the sake of synchronization between motors, has been growing recently. In this example, the consensus of a multi-motor system velocity is studied, which is inevitable to avoid damage to the product in industrial applications.

A multi-motor system containing three PMSMs with equivalent circuit, depicted in Figure 1, is considered. Due to the fractional behavior of capacitor and inductor, The fractional-order model of PMSM with parameter definitions in Table 1, is given as follows in [54]

$$\begin{cases} U_d - E = R(I_d + T_l\, D^{\zeta} I_d) \\ I_d - I_{dl} = \frac{T_m}{R} D^{\eta} E \end{cases} . \tag{65}$$

Transfer function between voltage and current is achieved by Laplace transform on either side of first equation of *(65)*, as follows

$$\frac{I_d(s)}{U_d(s)-E(s)} = \frac{1/R}{T_l s^{\zeta}+1}. \tag{66}$$

then, the Laplace transform on second equation of *(65)* results the transfer function between the current and electromotive force as,

$$\frac{E(s)}{I_d(s)-I_{dl}(s)} = \frac{R}{T_m s^{\eta}}, \tag{67}$$

The block diagram of PMSM is shown in Figure 2 considering that $v = E/C_e$ is the motor velocity.

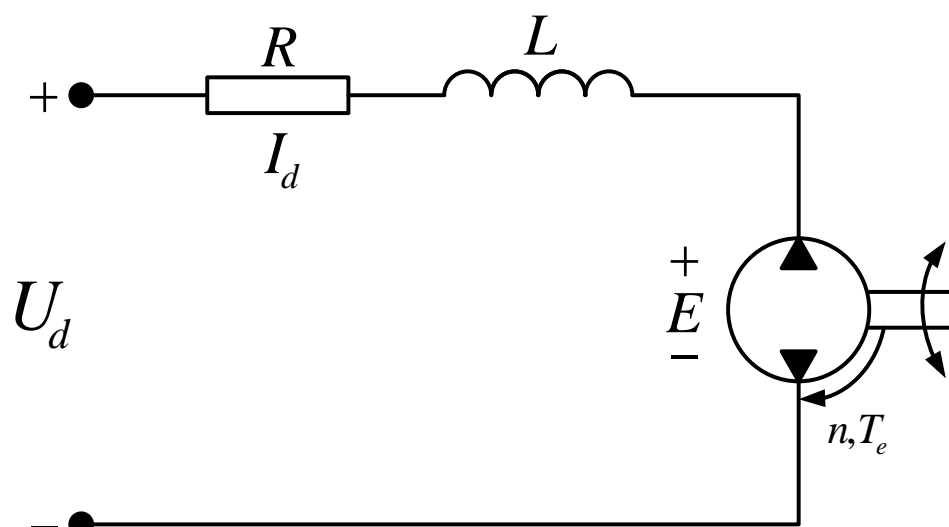

**Figure 1.** Equivalent circuit of synchronous motor

Finally, the transfer function of PMSM velocity control can be obtained as follows

$$G_{\zeta,\varsigma}(s) = \frac{1/C_e}{T_m T_l s^{\zeta+\eta}+T_m s^{\zeta}+1}. \tag{68}$$

**Table 1.** PMSM model terminology.

| Parameter | Description | Parameter | Description |
|---|---|---|---|
| $U_d$ | Armature voltage | $T_m$ | $T_m = GD^2R/375C_eC_m$ |
| $E$ | Back electromotive force | $C_e$ | Electromotive force coefficient of motor |
| $R$ | Stator resistance | $C_m$ | Torque constant ($C_m = (30/\pi)C_e$) |
| $I_d$ | Armature current | $\zeta$ | Fractional-order of inductor |
| $T_l$ | $L/R$ | $\eta$ | Fractional-order of capacitor |
| $I_{dl}$ | External load current | | |

According to the identified model of PMSM presented in [55]

$$G(s) = \frac{6.28}{0.00078s^{1.74}+0.097s^{0.87}+1}\ , \tag{69}$$

the pseudo-state space representation of the uncertain model *(68)* with additional nonlinear term is as follows

$$\begin{aligned} &D^{0.87}x_i(t) = \left(\tilde{A}+\Delta\tilde{A}\right)x_i(t) + B_iU_{di}(t) + \boldsymbol{\phi}(x_i(t), U_{di}),\\ &y_i(t) = \tilde{C}x_i(t), \quad i = 1,2,3, \end{aligned} \tag{70}$$

where $x_i = [v \quad D^{0.87}v]^T$, and

$$\begin{aligned} &\tilde{A} = \begin{bmatrix} 0 & 1 \\ -1282 & -124.3 \end{bmatrix},\ B_i = \begin{bmatrix} 0 \\ 6.28 \end{bmatrix},\ \tilde{C} = [1 \quad 0]\\ &\boldsymbol{\phi}(x_i(t), U_{di}) = \begin{bmatrix} sin\left(x_{i2}U_{di}(t)\right) + 0.5\,sin(x_{i2}) \\ -\,sin(x_{i1}) \end{bmatrix}. \end{aligned} \tag{71}$$

The uncertainty parameters are considered as

$$\widetilde{M} = \begin{bmatrix} 0.3 \\ 0 \end{bmatrix},\quad \tilde{R} = [-0.6 \quad 0.5],\quad J = 1,\quad \delta_i = sin\, i\pi/6, \tag{72}$$

where $\delta = \sin\pi/2 = 1$.

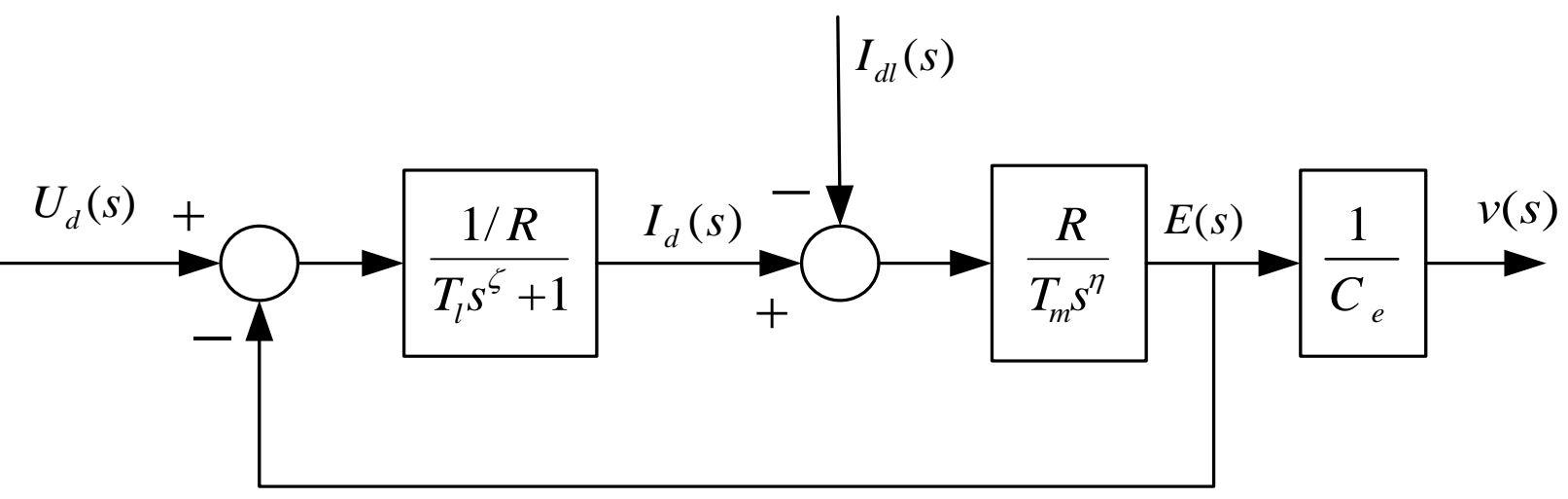

**Figure 2.** Block diagram of PMSM

The topology structure of the multi-agent system is demonstrated in Figure 3 with corresponding Laplacian matrix $L$, obtained according to *(14)*.

$$L = \begin{bmatrix} 1 & -1 & 0 \\ 0 & 1 & -1 \\ -1 & 0 & 1 \end{bmatrix}.$$

The controller uncertainty parameters are given in Table 2. To solve the consensus problem of uncertain system *(70)*, the dynamic and static controllers are designed using Theorem 2 . The resulted controllers are tabulated in Table 3. State trajectories of system resulted by the controllers of Table 3, for $n_c = 0,1,2$ are illustrated in Figure 4. According to the figures, motors' rotational velocity reaches consensus asymptotically in multi-motor system. Furthermore, as the controller is assumed to be dynamic, solver has more degree of freedom to find a feasible solution of inequality, because of additional parameters [45]. Although increasing controller order causes a slight improvement in system response, low order controller has desirable performance as well.

**Table 2.** Controller uncertainty parameters

| | **Order of controller** | $\boldsymbol{F(t)}$ | $\boldsymbol{D}$ | $\boldsymbol{E}$ |
|---|---|---|---|---|
| $\mathcal{A}_c$ | **1** | $\boldsymbol{diag}(sin(t), 0.5\,cos(3t), -cos(t))$ | $\begin{bmatrix} 0.2 & & \mathbf{0} \\ & -0.4 & \\ \mathbf{0} & & 0.9 \end{bmatrix}$ | $\begin{bmatrix} 0.4 & & \mathbf{0} \\ & 0.9 & \\ \mathbf{0} & & -0.1 \end{bmatrix}$ |
| | **2** | $\boldsymbol{diag}(sin(0.1t), cos(5t), sin(0.3t), sin(0.1t), cos(0.4t), sin(0.1t))$ | $\begin{bmatrix} 3.5 & 1.9 & & & & \\ 4.6 & 8.6 & & & & \mathbf{0} \\ & & 8.9 & 7.8 & & \\ & & 9.6 & 8.1 & & \\ & \mathbf{0} & & & 0 & 2.4 \\ & & & & 5.1 & 0 \end{bmatrix}$ | $\begin{bmatrix} 3.8 & 5.1 & & & & \\ 0.8 & 5.2 & & & & \mathbf{0} \\ & & 0.8 & 0.5 & & \\ & & 0.4 & 0.5 & & \\ & & & & 0.9 & 0 \\ & & & & 0.1 & 0.2 \end{bmatrix}$ |
| $\mathcal{B}_c$ | **1** | $\boldsymbol{diag}(cos(3t), sin(t), 0.2\,cos(t))$ | $\begin{bmatrix} 0.2 & & \mathbf{0} \\ & 0.4 & \\ \mathbf{0} & & -0.5 \end{bmatrix}$ | $\begin{bmatrix} 0.7 & & \mathbf{0} \\ & 0 & \\ \mathbf{0} & & -0.2 \end{bmatrix}$ |
| | **2** | $\boldsymbol{diag}(sin(t), cos(t), sin(t))$ | $\begin{bmatrix} 9 & & \\ 10 & & \mathbf{0} \\ & 2.4 & \\ & 9.5 & \\ \mathbf{0} & & 7.1 \\ & & 9.3 \end{bmatrix}$ | $\begin{bmatrix} 0.2 & & \mathbf{0} \\ & 0.6 & \\ \mathbf{0} & & 0.5 \end{bmatrix}$ |
| $\mathcal{C}_c$ | **1** | $\boldsymbol{diag}(sin(0.5t), cos(t), cos\,(t))$ | $\begin{bmatrix} 0.4 & & \mathbf{0} \\ & -0.1 & \\ \mathbf{0} & & 0.8 \end{bmatrix}$ | $\begin{bmatrix} 0.3 & & \mathbf{0} \\ & -0.4 & \\ \mathbf{0} & & 0.5 \end{bmatrix}$ |
| | **2** | $\boldsymbol{diag}(sin(t), cos(t), sin(t))$ | $\begin{bmatrix} 5.8 & & \mathbf{0} \\ & 8.8 & \\ \mathbf{0} & & 6 \end{bmatrix}$ | $\begin{bmatrix} 3.2 & & \\ 5.1 & & \mathbf{0} \\ & 0.4 & \\ & 0 & \\ \mathbf{0} & & -0.4 \\ & & -0.5 \end{bmatrix}^T$ |
| $\mathcal{D}_c$ | **1** | $\boldsymbol{diag}(cos(2t), sin(t), sin(t))$ | $\begin{bmatrix} 0.7 & & \mathbf{0} \\ & 0.7 & \\ \mathbf{0} & & -0.1 \end{bmatrix}$ | $\begin{bmatrix} 0.4 & & \mathbf{0} \\ & 0.4 & \\ \mathbf{0} & & 0.6 \end{bmatrix}$ |
| | **2** | $\boldsymbol{diag}(sin(0.5t), cos(10t), sin(0.1t))$ | $\begin{bmatrix} 7.2 & & \mathbf{0} \\ & 6.4 & \\ \mathbf{0} & & 2.8 \end{bmatrix}$ | $\begin{bmatrix} 1.3 & & \mathbf{0} \\ & 0.7 & \\ \mathbf{0} & & 3.2 \end{bmatrix}$ |

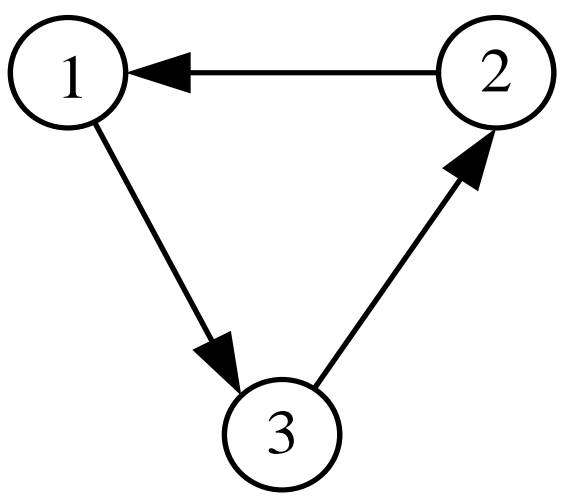

**Figure 3.** Network topology of PMSM.

Besides, control effort trajectories of proposed controllers in Table 3, are plotted in Figure 5. Perusing these figures, we can conclude that increasing controller order, leads to a perceptible enhancement in consensus rate. Moreover, to have a meaningful comparison we have investigated the results with consensus time $CT_p$ which refers to a time that for $t \geq CT_p$, $\$\sum_{i=1}^{N} e_i(t) \leq p \times 10^{-2} \sum_{i=1}^{N} e_{i0}$, where $e_i(t)$ is the consensus error defined in Definition 2 and $e_{i0} = e_i(0)$. Time $CT_p$ is when consensus is achieved within the $p\%$ of the initial consensus error [56]. As it is shown in figures the $CT_2$ for the resulted dynamic controllers with the order 2, 1 from theorem 2 and the static controller [20], [21] are obtained as $0.11$, $0.18$ and $0.13$. Comparing $CT_2$ shows that the second order dynamic controller has the better outcome in comparison with the first order and static controllers. Analyzing the $CT_2$ and controller effort for the first order and static controllers demonstrates a competitive result. The static controller illustrates a preferable outcome from the $\mathrm{CT}_2$ viewpoint, on the other hand it has a greater control effort. Simply, static controller achieves a better $CT_2$ by consuming more energy. But the first order controller with a normal control effort has a desired $CT_2$ as well.

**Table 3.** Controller parameters obtained by Theorem 2 .

| **Order of controller** | | **Agent 1** | **Agent 2** | **Agent 3** |
|---|---|---|---|---|
| **0** | $\boldsymbol{D_c}$ | $-105.33$ | $-59.70$ | $-59.95$ |
| **1** | $\boldsymbol{A_c}$ | $-57.45$ | $-66.14$ | $-22.43$ |
| | $\boldsymbol{B_c}$ | $16.28$ | $3.54$ | $11.03$ |
| | $\boldsymbol{C_c}$ | $-15.26$ | $-7.45$ | $-3.14$ |
| | $\boldsymbol{D_c}$ | $-83.74$ | $-45.31$ | $-45.82$ |
| **2** | $\boldsymbol{A_c}$ | $\begin{bmatrix} -50.58 & -34.24 \\ -34.29 & -78.37 \end{bmatrix}$ | $\begin{bmatrix} -61.71 & -42.54 \\ -42.17 & -84.95 \end{bmatrix}$ | $\begin{bmatrix} -42.05 & -23.86 \\ -23.85 & -62.63 \end{bmatrix}$ |
| | $\boldsymbol{B_c}$ | $\begin{bmatrix} 13.33 \\ 16.95 \end{bmatrix}$ | $\begin{bmatrix} 9.33 \\ 3.03 \end{bmatrix}$ | $\begin{bmatrix} 15.34 \\ 10.30 \end{bmatrix}$ |
| | $\boldsymbol{C_c}$ | $[-50.75 \quad -18.36]$ | $[-16.74 \quad -40.21]$ | $[-38.32 \quad -11.84]$ |
| | $\boldsymbol{D_c}$ | $-74.25$ | $-53.82$ | $-51.79$ |

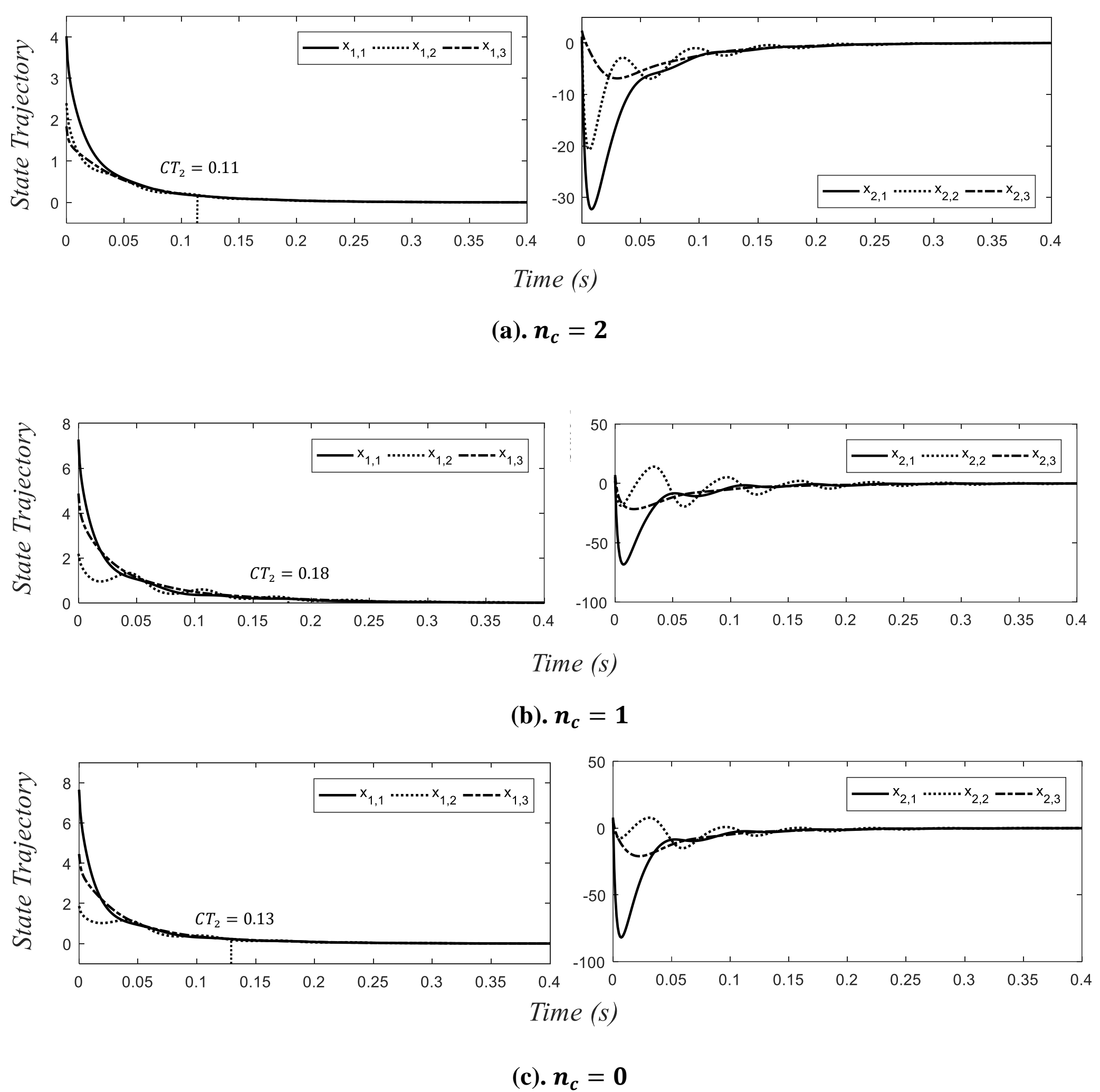

**(a). $\boldsymbol{n_c = 2}$**

**(b). $\boldsymbol{n_c = 1}$**

**(c). $\boldsymbol{n_c = 0}$**

**Figure 4.** Velocity of PMSM using the proposed dynamic output feedback controller in Theorem 2 .

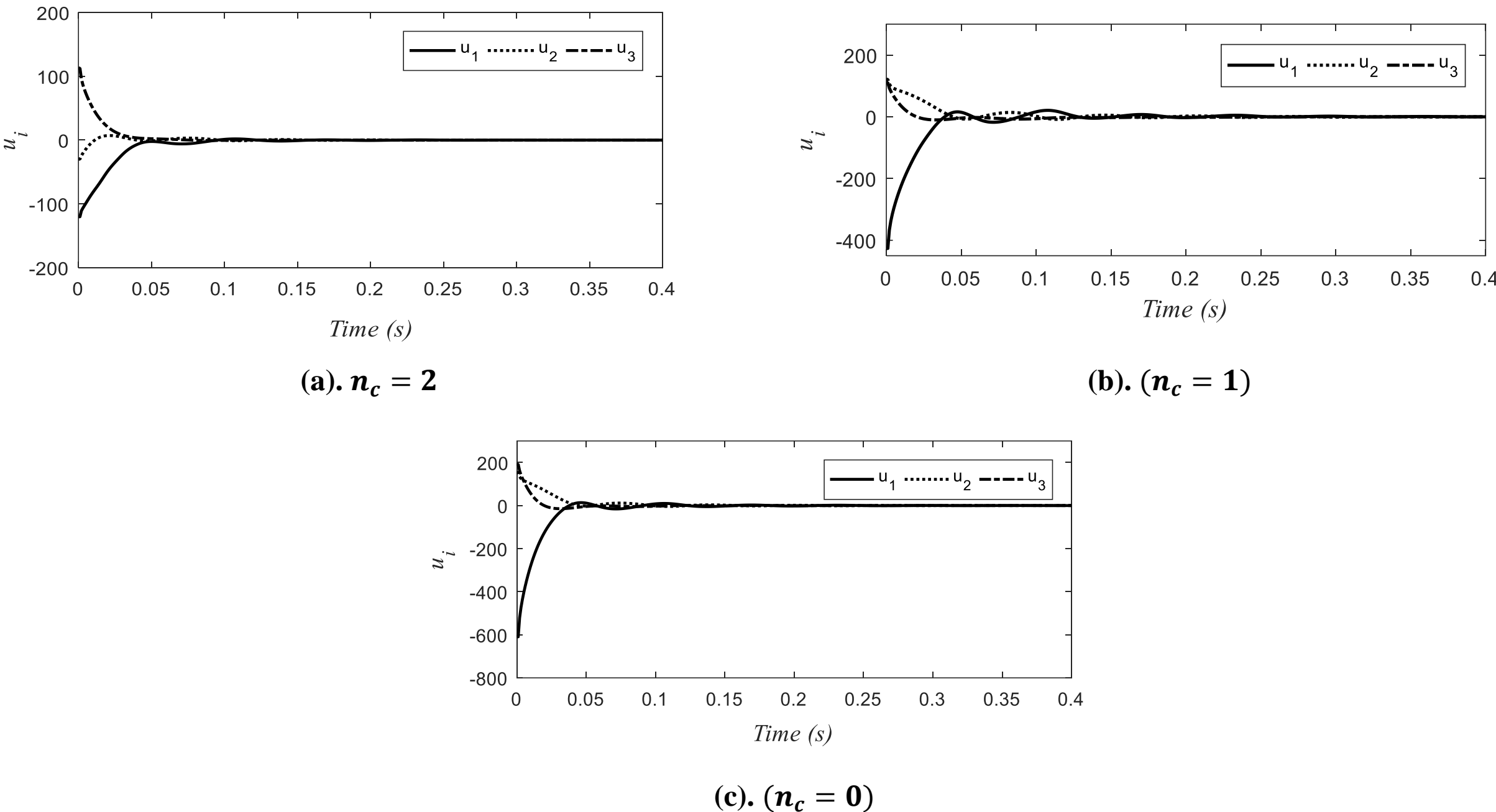

**(a).** $\boldsymbol{n_c = 2}$ **(b).** $(\boldsymbol{n_c = 1})$

**(c).** $(\boldsymbol{n_c = 0})$

**Figure 5.** Control effort of proposed controller in Theorem 2 .

### B. Numerical example

A connected network with four agents is considered as shown in Figure 6. The dynamics of each agent are represented as follows

$$\begin{aligned} & D^{\alpha}x_i(t) = \left(\tilde{A} + \Delta\tilde{A}_i\right)x_i(t) + B_i U_{di}(t) \ + \boldsymbol{\phi}(x_i(t), U_{di}),, \\ & y_i(t) = \tilde{C}x_i(t), \qquad\qquad\qquad i = 1, 2, 3, \end{aligned} \tag{73}$$

where

$$\tilde{A} = \begin{bmatrix} -1 & 1 & 0 & 0 \\ 1 & -3 & 0 & 1 \\ 0 & 0 & 0 & 1 \\ 0 & 0 & -1 & 0 \end{bmatrix}, \ B_i = \begin{bmatrix} 1 \\ 1 \\ 0 \\ 1 \end{bmatrix}, \ \tilde{C} = \begin{bmatrix} 1 \\ 0 \\ 1 \\ 0 \end{bmatrix}^T,$$

$$\boldsymbol{\phi}(x_i(t), U_{di}) = \begin{bmatrix} sin\left(x_{i2}U_{di}(t)\right) + 0.5\,sin(x_{i2}) \\ -\,sin(x_{i1}) \end{bmatrix}, \tag{74}$$

$$\widetilde{M} = [0.2 \quad 0 \quad -0.1 \quad 0.3]^T, \ \ \tilde{R} = [0 \quad 0.2 \quad 0.4 \quad -0.2], \quad J = 1,$$

$$\delta_1 = 0.5, \ \delta_2 = -0.4, \ \delta_3 = 0.1, \ \delta_4 = 0.8,$$

where $\delta = 0.8$. The corresponding Laplacian matrix of the mentioned network is

$$L = \begin{bmatrix} 2 & -1 & 0 & -1 \\ -1 & 2 & -1 & 0 \\ 0 & -1 & 2 & -1 \\ -1 & 0 & -1 & 2 \end{bmatrix}.$$

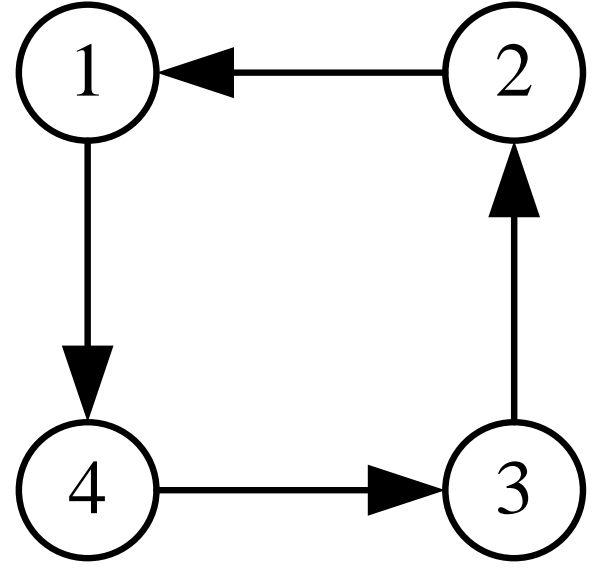

**Figure 6.** Network topology o FOMAS.

This example is solved for $\alpha = 0.8$, with controller uncertainty parameters, given in Table 4, where, the resulted controller parameters obtained by Theorem 2 and Corollary 1 for system (73), with and without the nonlinear term respectively, are presented in Table 5. Figure 7 and Figure 8 illustrate the state trajectories of agents for initial states using the output feedback controllers presented in Theorem 2 and Corollary 1 , besides control effort of proposed controllers depicted in Figure 9 and Figure 10. Obviously, the corresponding trajectories of all agents asymptotically reach an agreement.

**Table 4.** Controller uncertainty parameters of system *(73)*

| | **Order of controller** | $\boldsymbol{F(t)}$ | $\boldsymbol{D}$ | $\boldsymbol{E}$ |
|---|---|---|---|---|
| $\mathcal{A}_c$ | **1** | $\boldsymbol{diag}(sin(t), cos(5t), sin(3t), sin(t))$ | $\begin{bmatrix} 1.25 & & & \mathbf{0} \\ & 0.25 & & \\ & & 0.4 & \\ \mathbf{0} & & & 1.75 \end{bmatrix}$ | $\begin{bmatrix} 0.8 & & & \mathbf{0} \\ & 0.7 & & \\ & & 0.1 & \\ \mathbf{0} & & & 0.6 \end{bmatrix}$ |
| | **2** | $\boldsymbol{diag}(sin(t), cos(5t), sin(3t), sin(t), cos(4t), sin(t), sin(t), sin\ (t))$ | $\begin{bmatrix} 1.25 & 1 & & & & & & \mathbf{0} \\ 1.25 & 2.5 & & & & & & \\ & & 0.25 & 1.5 & & & & \\ & & 1.25 & 1.5 & & & & \\ & & & & 0.4 & 0 & & \\ & & & & 2 & 1.5 & & \\ & & & & & & 1.75 & 1.5 \\ \mathbf{0} & & & & & & 2.25 & 1.75 \end{bmatrix}$ | $\begin{bmatrix} 0.8 & 0.6 & & & & & & \mathbf{0} \\ 0.5 & 0.6 & & & & & & \\ & & 0.7 & 0.2 & & & & \\ & & 0.3 & 0.3 & & & & \\ & & & & 0 & 0.9 & & \\ & & & & 0.1 & 0.8 & & \\ & & & & & & 0.6 & 4 \\ \mathbf{0} & & & & & & 5.7 & 6.1 \end{bmatrix}$ |
| $\mathcal{B}_c$ | **1** | $\boldsymbol{diag}(sin(t), cos(t), sin(t), sin(t))$ | $\begin{bmatrix} 2 & & & \mathbf{0} \\ & 0.75 & & \\ & & 2.25 & \\ \mathbf{0} & & & 0.75 \end{bmatrix}$ | $\begin{bmatrix} 0.5 & & & \mathbf{0} \\ & 0.5 & & \\ & & 0.7 & \\ \mathbf{0} & & & 0.4 \end{bmatrix}$ |
| | **2** | $\boldsymbol{diag}(sin(t), cos(t), sin(t), sin(t))$ | $\begin{bmatrix} 2 & & & \mathbf{0} \\ 0.5 & & & \\ & 0.75 & & \\ & 0.5 & & \\ & & 2 & \\ & & 2.25 & \\ & & & 0.75 \\ \mathbf{0} & & & 0.15 \end{bmatrix}$ | $\begin{bmatrix} 0.5 & & & \mathbf{0} \\ & 0.5 & & \\ & & 0.7 & \\ \mathbf{0} & & & 0.4 \end{bmatrix}$ |
| $\mathcal{C}_c$ | **1** | $\boldsymbol{diag}(sin(t), cos(t), sin(t), sin(t))$ | $\begin{bmatrix} 0.5 & & & \mathbf{0} \\ & 1.75 & & \\ & & 1.75 & \\ \mathbf{0} & & & 0 \end{bmatrix}$ | $\begin{bmatrix} 25 & & & \mathbf{0} \\ & 10 & & \\ & & 1.8 & \\ \mathbf{0} & & & 2.4 \end{bmatrix}$ |
| | **2** | $\boldsymbol{diag}(sin(t), cos(t), sin(t), sin(t))$ | $\begin{bmatrix} 0.5 & & & \mathbf{0} \\ & 1.75 & & \\ & & 1.75 & \\ \mathbf{0} & & & 0 \end{bmatrix}$ | $\begin{bmatrix} 40 & & & \mathbf{0} \\ 25 & & & \\ & 10 & & \\ & 5 & & \\ & & 1.8 & \\ & & 4.6 & \\ & & & 2.4 \\ \mathbf{0} & & & 5 \end{bmatrix}^T$ |
| $\mathcal{D}_c$ | **1** | $\boldsymbol{diag}(sin(t), cos(t), sin(5t), sin(t))$ | $\begin{bmatrix} 1 & & & \mathbf{0} \\ & 0 & & \\ & & 0.75 & \\ \mathbf{0} & & & 0.75 \end{bmatrix}$ | $\begin{bmatrix} 2.1 & & & \mathbf{0} \\ & -0.7 & & \\ & & 1.4 & \\ \mathbf{0} & & & 0.2 \end{bmatrix}$ |
| | **2** | $\boldsymbol{diag}(sin(t), cos(t), sin(5t), sin(t))$ | $\begin{bmatrix} 1 & & & \mathbf{0} \\ & 0 & & \\ & & 0.75 & \\ \mathbf{0} & & & 0.75 \end{bmatrix}$ | $\begin{bmatrix} 2.1 & & & \mathbf{0} \\ & -0.7 & & \\ & & 1.4 & \\ \mathbf{0} & & & 0.2 \end{bmatrix}$ |

The consensus error (defined in (17)) indices of proposed static and dynamic controllers are summarized in Table 6 where the results indicate that increasing controller order reduces the settling time of the controller effort. Since the controller objective is the consensus of the corresponding states of agents, the vanishing of the consensus error, indicates that the consensus of agents is achieved. Thus, there is a direct relationship between the controller order and the rate of the consensus. Besides increasing the controller, order leads to a slight decrement in consensus error indices, and a quicker and more efficient consensus.

**Table 5.** Controller parameters obtained by (I) Theorem 2 and, (II) Corollary 1

| | Order of controller | | Agent 1 | Agent 2 | Agent 3 | Agent 4 |
|---|---|---|---|---|---|---|
| I | 1 | $\boldsymbol{A_c}$ | 56.9751 | 50.2846 | 31.2825 | 25.4864 |
| | | $\boldsymbol{B_c}$ | −3.3796 | 3.1756 | −5.4820 | −2.5225 |
| | | $\boldsymbol{C_c}$ | 1.3433 | −0.2796 | 1.5278 | −5.0534 |
| | | $\boldsymbol{D_c}$ | −45.9537 | −50.5748 | −44.8136 | −42.0247 |
| | 2 | $\boldsymbol{A_c}$ | $\begin{bmatrix} -32.0576 & -14.8811 \\ -14.9358 & -36.1773 \end{bmatrix}$ | $\begin{bmatrix} -17.8004 & -13.9017 \\ -13.8483 & -24.0191 \end{bmatrix}$ | $\begin{bmatrix} -24.5997 & -21.2761 \\ -22.9911 & -42.9351 \end{bmatrix}$ | $\begin{bmatrix} -24.5167 & -25.4635 \\ -24.4884 & -36.1116 \end{bmatrix}$ |
| | | $\boldsymbol{B_c}$ | $[1.537 \quad 1.624]^T$ | $[-0.679 \quad -0.792]^T$ | $[-0.525 \quad -0.694]^T$ | $[-1.988 \quad -2.366]^T$ |
| | | $\boldsymbol{C_c}$ | $[-2.1540 \quad -4.7027]$ | $[6.5429 \quad -7.9858]$ | $[-6.0084 \quad -4.1185]$ | $[3.3085 \quad 3.0719]$ |
| | | $\boldsymbol{D_c}$ | −27.9259 | −27.5037 | −27.1510 | −8.0702 |
| II | 2 | $\boldsymbol{A_c}$ | $\begin{bmatrix} -38.45 & -0.15 \\ -0.15 & -38.45 \end{bmatrix}$ | $\begin{bmatrix} -36.97 & -0.2 \\ -0.23 & -35.21 \end{bmatrix}$ | $\begin{bmatrix} -33.12 & -1.54 \\ -1.54 & -40.45 \end{bmatrix}$ | $\begin{bmatrix} -29.88 & -0.65 \\ -0.73 & -27.03 \end{bmatrix}$ |
| | | $\boldsymbol{B_c}$ | $[4.67 \quad 4.67]^T$ | $[4.67 \quad 4.67]^T$ | $[4.67 \quad 4.67]^T$ | $[4.67 \quad 4.67]^T$ |
| | | $\boldsymbol{C_c}$ | $[0.04 \quad 0.11]$ | $[0.14 \quad 0.14]$ | $[0.07 \quad 0.07]$ | $[0.16 \quad 0.16]$ |
| | | $\boldsymbol{D_c}$ | −3.75 | −3.75 | −3.75 | −3.75 |

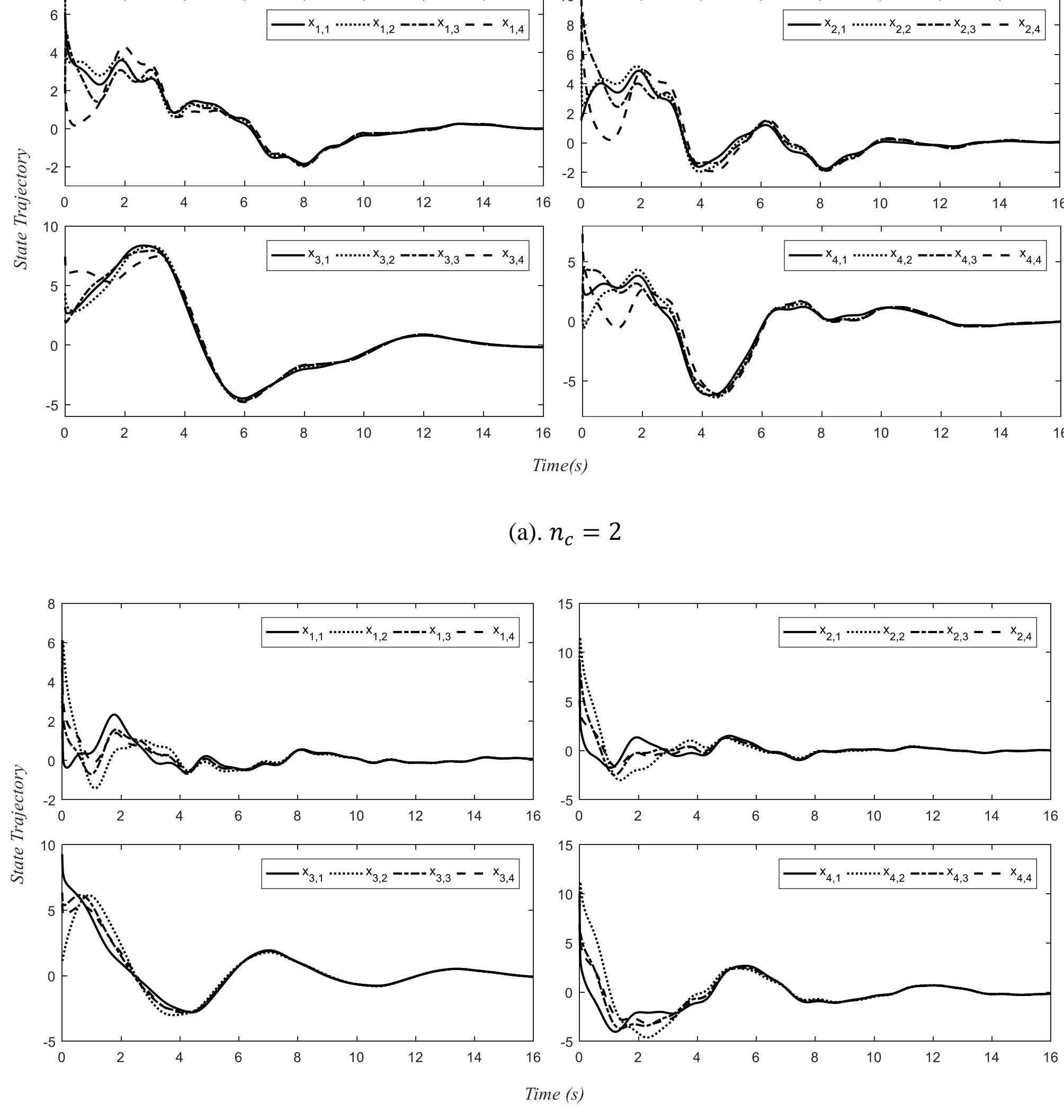

(a). $n_c = 2$

(b). $n_c = 1$

**Figure 7.** State trajectories of multi-agent system shown in Figure 6 using the proposed controller in Theorem 2 **.**

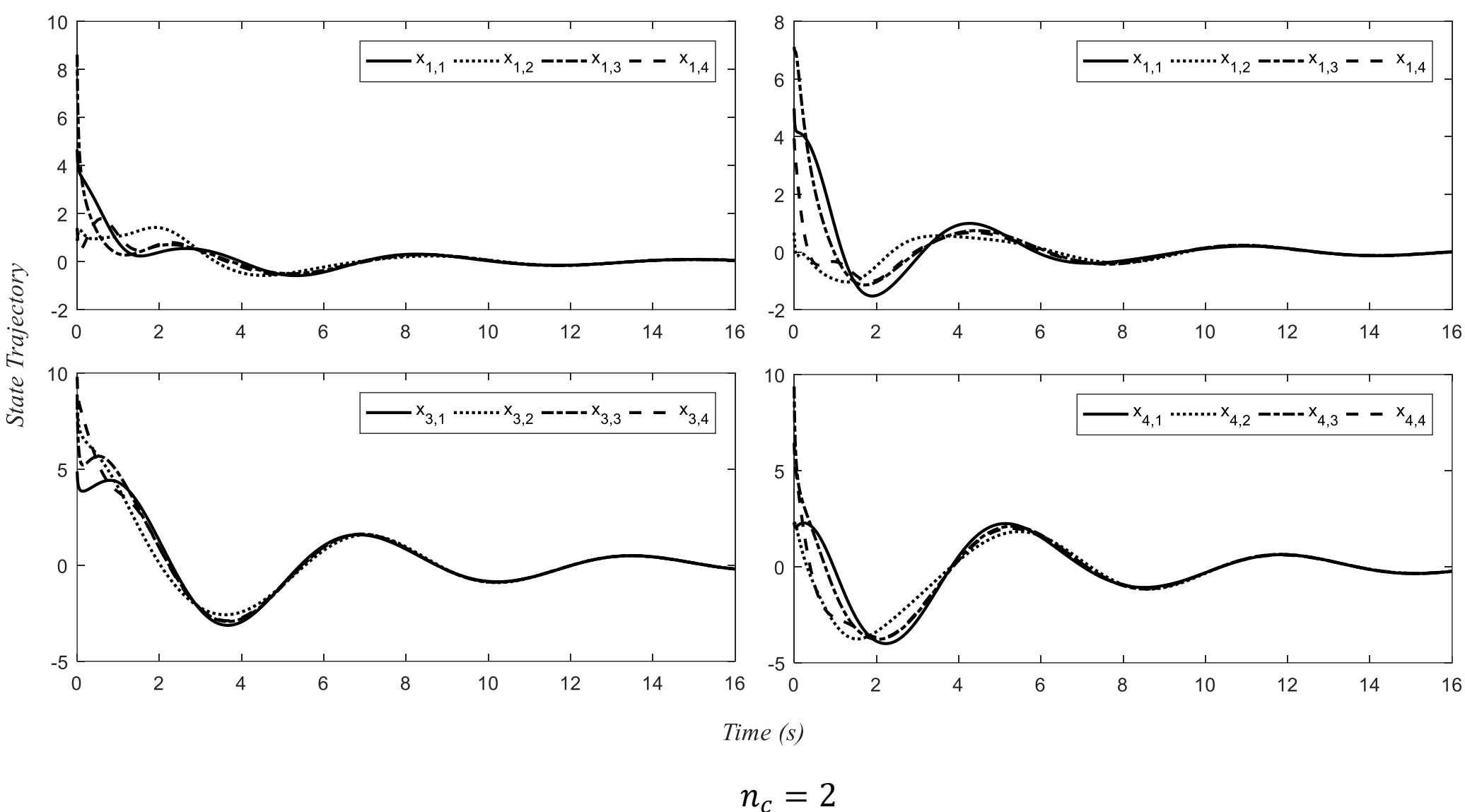

$n_c = 2$

**Figure 8.** State trajectories of multi-agent system shown in Figure 6 using the proposed controller in Corollary 1 **.**

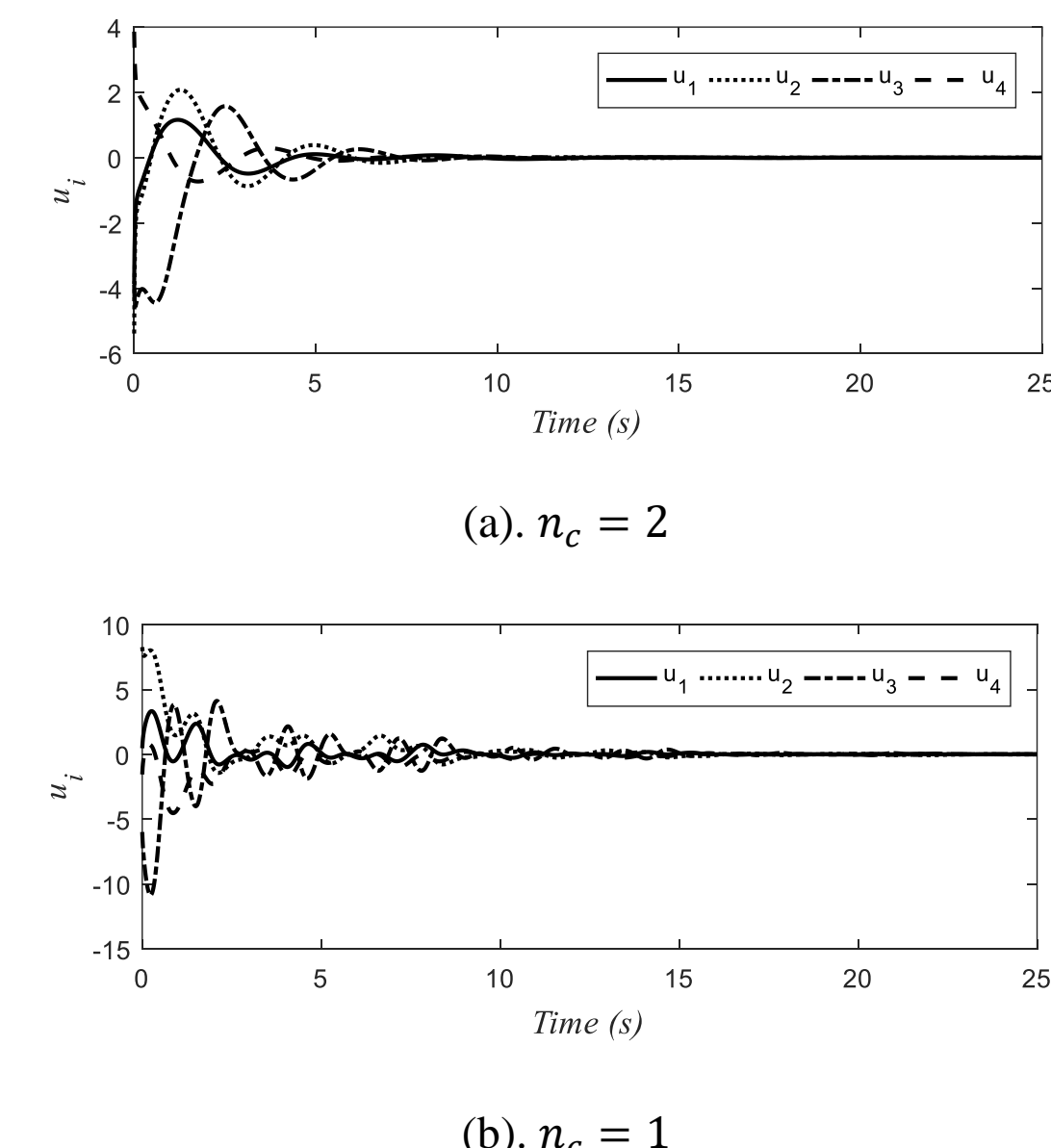

(a). $n_c = 2$

(b). $n_c = 1$

**Figure 9** Control effort of proposed controller in Theorem 2 .

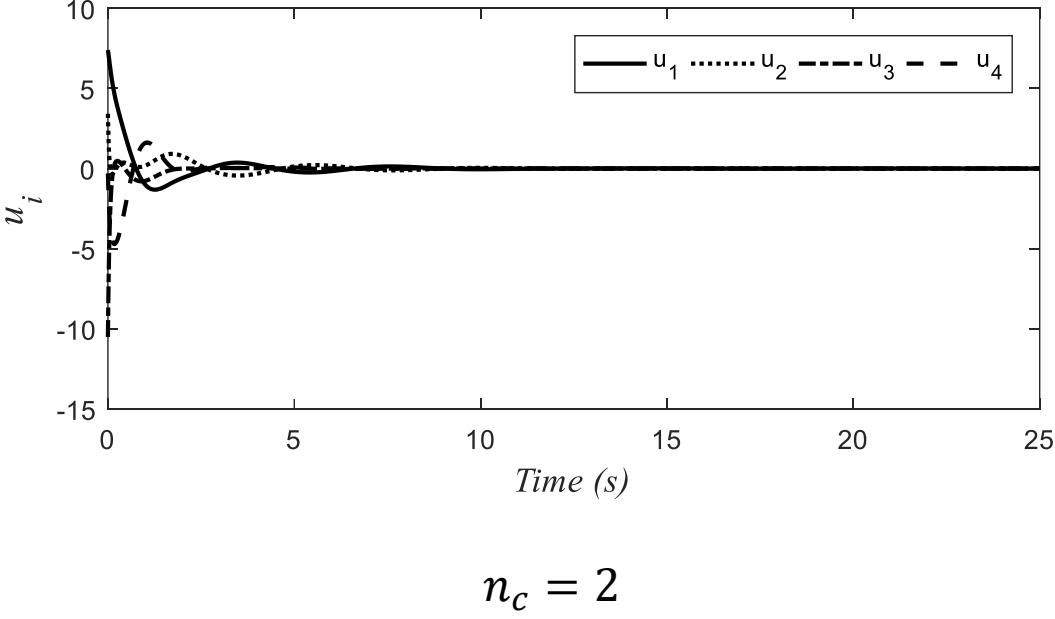

$n_c = 2$

**Figure 10** Control effort of proposed controller in Corollary 1 .

**Table 6.** Consensus error indices for proposed method in (I) Theorem 2 and (II) Corollary 1

| **Controller** | **order** | | **ISE** | **IAE** | **ITSE** | **ITAE** |
|---|---|---|---|---|---|---|
| **I** | **1** | Agent 1 | 1878 | 88.31 | 2087 | 183.98 |
| | | Agent 2 | 1358 | 72.32 | 1297 | 141.81 |
| | | Agent 3 | 714.24 | 54.14 | 742.71 | 111.09 |
| | | Agent 4 | 644.30 | 51.02 | 650.05 | 104.30 |
| | **2** | Agent 1 | 424.96 | 46.31 | 448.87 | 146.36 |
| | | Agent 2 | 1841 | 82.54 | 1638 | 171.93 |
| | | Agent 3 | 488.49 | 46.74 | 485.72 | 123.43 |
| | | Agent 4 | 558.78 | 47.61 | 506.28 | 113.43 |
| **II** | **2** | Agent 1 | 592.83 | 45.22 | 430.91 | 97.57 |
| | | Agent 2 | 570.24 | 47.79 | 501.49 | 112.73 |
| | | Agent 3 | 468.50 | 35.23 | 232.63 | 71.47 |
| | | Agent 4 | 392.48 | 33.55 | 241.12 | 71.07 |

## V. CONCLUSION

In this paper, decentralized static and dynamic output feedback non-fragile controllers for the consensus of a class of nonlinear fractional-order multi-agent systems with positive real uncertainty are proposed. First, a new FOMAS with transformed states is defined, in which the stability of this new system is equivalent to the consensus of the main system. Second, sufficient conditions for the stability of new system using fractional-order systems stability theorems and Schur complement are obtained in the form of linear matrix inequalities. Third, the controller unknown parameters obtained by solving matrix inequalities. Eventually, some numerical examples are presented to illustrate the effectiveness of the proposed dynamic output feedback controller design methods for FOMASs.